\documentclass{article}
\usepackage[T2A]{fontenc}
\usepackage[cp1251]{inputenc}
\usepackage[english,russian]{babel}
\usepackage[tbtags]{amsmath}
\usepackage{amsfonts,amssymb,mathrsfs,amscd}
\numberwithin{equation}{section}
\newtheorem{remark}{Remark}[section]
\newtheorem{lemma}{Lemma}[section]
\newtheorem{theorem}{Theorem}[section]
\newtheorem{definition}{Definition}[section]
\newtheorem{corollary}{Corollary}[section]

\newcommand{\sign}{\mathop{\rm sign}}

\renewcommand{\span}{\mathop{\rm span}}
\renewcommand{\Im}{\mathop{\rm Im}}
\renewcommand{\Re}{\mathop{\rm Re}}
\newcommand{\Ker}{\mathop{\rm Ker}}
\newcommand{\supp}{\mathop{\rm supp}}
\newcommand{\rank}{\mathop{\rm rank}}
\newcommand{\col}{\mathop{\rm col}}

\begin{document}
\begin{Large}
\thispagestyle{empty}
\begin{center}
{\bf Characteristic function of pencils. Model representations of a quadratic operator pencil\\
\vspace{5mm}
V. A. Zolotarev}\\

B. Verkin Institute for Low Temperature Physics and Engineering
of the National Academy of Sciences of Ukraine\\
47 Nauky Ave., Kharkiv, 61103, Ukraine

Department of Higher Mathematics and Informatics, V. N. Karazin Kharkov National University \\
4 Svobody Sq, Kharkov, 61077,  Ukraine
\vspace{5mm}

{\it Dedicated to my teacher, M. S. Liv$\breve{s}$ic}

\end{center}
\vspace{5mm}

{\small {\bf Abstract.} Notion of an open system of second order is introduced. Characteristic function for such an open system is obtained. Model representations of a quadratic non-self-adjoint operator pencil are found.}
\vspace{5mm}

%{\it Mathematics
%Subject Classification 2020:} 34L10, 34L15.\\

\vspace{5mm}

\begin{center}
{\bf Introduction}
\end{center}
\vspace{5mm}

Open system corresponding to a non-self-adjoint bounded operator was introduced by M. S. Liv$\breve{\rm s}$ic. This open system is of first order and underlies two directions closely related to the theory of non-self-adjoint operators, namely: dilation theory (Nagy -- Foias) and scattering theory (Lax -- Phillips) \cite{24}. It is an open system that leads to characteristic function of a non-self-adjoint operator which is the main analytic instrument for obtaining triangular (M. S. Liv$\breve{\rm s}$ic) and functional (Nagy -- Foias) models.

As a rule, physical laws are described by second order equations, therefore it is natural to consider an open system defined by the second order equation
$$i\ddot{h}+B\dot{h}+Ah=f$$
where $B$ is a self-adjoint operator characterizing dissipation in the system and the operator $A$ is a non-self-adjoint bounded operator (e.g., of finite range of non-hermiticity). This equation can be reduced by a standard method to a first order equation and then the operator $B$ is included in the main operator as a separate block. As a result, information about dissipation in the system is lost. This fact is the main motivation for study of second order open systems.

The paper consists of four sections. The first section gives definition of an open system of second order and characteristic function
$$S_\Delta(\lambda,B)=I-i\varphi(\Lambda^2I+\lambda B+A)^{-1}\varphi^*\sigma$$
($A-A^*=i\varphi^*\sigma\varphi$). Its properties and operation of coupling of open systems are described.

Section 2 gives analogue of a well-known statement on representation of a self-adjoint operator by an operator of multiplication by independent variable.

As is seen from the form of characteristic function $S_\Delta(\lambda,B)$, it is necessary to study spectral properties of a pencil $L(\lambda)=\lambda^2I+\lambda B+A$. According to M. G. Krein, to do this one has to factorize a pencil
$$L(\lambda)=(\lambda I-Y)(\lambda I-X)$$
and study its ``roots'' $X$, $Y$.

Section 3 constructs model representations of the operators $X$ and $Y$. It is shown that if one of the roots is the operator of multiplication by independent variable, the the second operator is generalization of Hilbert transform (or Stieltjes transform).

Section 4 calculates characteristic function $S_\Delta(\lambda,B)$ for models constructed in Section 3. Notice that Riemann boundary value problem is used in this calculation.

Unsolved problems in this subject are stated in the conclusion.

\section{Second order open systems}\label{s1}

{\bf 1.1.} Remind \cite{23, 24} that totality of Hilbert spaces $H$, $E$ and linear bounded operators $A$: $H\rightarrow H$, $\varphi$: $H\rightarrow E$, $\sigma$: $E\rightarrow E$ ($\sigma=\sigma^*$) is said to be a {\bf local colligation}
\begin{equation}
\Delta=(A,H,\varphi,E,\sigma)\label{eq1.1}
\end{equation}
if
\begin{equation}
A-A^*=i\varphi^*\sigma\varphi.\label{eq1.2}
\end{equation}
The open system $\mathcal{F}_\Delta=\{\mathcal{R}_\Delta,S_\Delta\}$ is commonly associated with a colligation $\Delta$. This system is given by the first order differential equation
\begin{equation}
\begin{array}{lll}
\displaystyle{\mathcal{R}_\Delta:\left\{
\begin{array}{lll}
i\dot{h}+Ah=\varphi^*\sigma u;\\
\left.h\right|_0=h_0;
\end{array}\right.}\\
S_\Delta:v=u-i\varphi h,
\end{array}\label{eq1.3}
\end{equation}
here $h=h(t)$ and $u=u(t)$, $v=v(t)$ are vector functions from $H$ and $E$ respectively; $\displaystyle{\dot{h}=\frac{dh}{dt}}$. For an open system \eqref{eq1.2}, the {\bf conservation law} is true,
\begin{equation}
\frac d{dt}\|h\|^2=\langle\sigma u,u\rangle-\langle\sigma v,v\rangle\label{eq1.4}
\end{equation}
on which construction of {\bf unitary dilation} \cite{23} of a semigroup of contractions $Z_t=\exp\{itA\}$ ($\sigma=I$, $t\in\mathbb{R}_+$) is based. Mapping $S_\Delta$  \eqref{eq1.3}, after the Fourier transform, becomes the {\bf characteristic function} $S_\Delta(\lambda)$,
\begin{equation}
S_\Delta(\lambda)=I-i\varphi(A-\lambda I)^{-1}\varphi^*\sigma,\label{eq1.5}
\end{equation}
of the colligation $\Delta$ \eqref{eq1.1}. Characteristic function was introduced by M. S. Livsic and is a main instrument \cite{23, 24, 35} of spectral analysis of a non-self-adjoint operator $A$.
\vspace{5mm}

{\bf 1.2.} Match colligation $\Delta$ \eqref{eq1.1} with a second order open system. Let $B$ be a linear bounded operator in $H$.

\begin{definition}\label{d1.1}
Pair of mappings $\mathcal{R}_\Delta(B)(h_0,h_1,u)=h$; $S_\Delta(h_0,h_1,u)=$ $(h_0(T),h_1(T),v)$ of the form
\begin{equation}
\mathcal{R}_\Delta(B):\left\{
\begin{array}{lll}
\ddot{h}+B\dot{h}+Ah=\varphi^*\sigma u;\\
\left.h\right|_0=h_0,\left.\dot{h}\right|_0=h_1;\quad(0\leq t\leq T<\infty);
\end{array}\right.\label{eq1.6}
\end{equation}
\begin{equation}
S_\Delta:v=u-i\varphi h\label{eq1.7}
\end{equation}
is said to be the {\bf second order open system} $\mathcal{F}_\Delta(B)=\{\mathcal{R}_\Delta(B),S_\Delta\}$ {\bf associated with a pair} $\{\Delta,B\}$. Here $h_0$, $h_1\in H$; $\displaystyle{\dot{h}=\frac d{dt}h}$, ${\displaystyle\ddot{h}=\frac{d^2h}{dt^2}}$; function $h$ in \eqref{eq1.7} is a solution to Cauchy problem \eqref{eq1.6}; $h_0(T)=h(T)$, $h_1(T)=\dot{h}(T)$.
\end{definition}

Operator $B$ in \eqref{eq1.6} describes dissipation (losses) of the process defined by equation \eqref{eq1.6}. Dissipation is caused by friction, viscosity, etc. In physics problems, operator $B$, as a rule, is self-adjoint, $B=B^*$, this is supposed hereinafter.

Equation \eqref{eq1.6}, due to doubling of initial space, $\widetilde{h}=\col\left[h,\dot{h}\right]$, is transformed into the first order equation $\dot{\widetilde{h}}+\widetilde{A}\widetilde{h}=\varphi^*\sigma u$, where the operators $A$ and $B$ are included in $\widetilde{A}$ as separate blocks. In such a representation operator $B$ characterizing impact of `environment' on the physical process is lost. This is reasonable enough motivation for studying the equations of \eqref{eq1.6} kind.

\begin{remark}\label{r1.1}
One can choose a transfer mapping of the form
\begin{equation}
S_\Delta:v=u-i\varphi(h+C\dot{h}),\label{eq1.7'}
\end{equation}
here $C$ is a linear bounded operator in $H$. Study of the systems \eqref{eq1.6}, \eqref{eq1.7'} is rather complicated and we confine ourselves to the case of $C=0$.
\end{remark}

\begin{theorem}
For the open system $\mathcal{F}_\Delta(B)$ \eqref{eq1.6}, \eqref{eq1.7} the {\bf conservation law} holds,
\begin{equation}
\langle\sigma u,u\rangle-\langle\sigma v,v\rangle=\frac d{dt}2\Im\langle\dot{h},h\rangle+2\Im\langle B\dot{h},h\rangle\label{eq1.8}
\end{equation}
where $B=B^*$.
\end{theorem}

P r o o f. Equations \eqref{eq1.2}, \eqref{eq1.6}, \eqref{eq1.7} imply
$$\langle\sigma u,u\rangle-\langle\sigma v,v\rangle=\langle\sigma u,u-v\rangle+\langle\sigma(u-v),u\rangle-\langle\sigma(u-v),(u-v)\rangle$$
$$=\langle\sigma u,i\varphi h\rangle+\langle i\varphi h,\sigma u\rangle-\langle\varphi^*\sigma\varphi h,h\rangle=-i\langle\ddot{h}+B\dot{h}+Ah,h\rangle$$
$$+i\langle h,\ddot{h}+B\dot{h}+Ah\rangle+i\langle(A-A^*)h,h\rangle=-i(\langle\ddot{h},h\rangle-\langle h,\ddot{h}\rangle)$$
$$-i(\langle B\dot{h},h\rangle-\langle h,B\dot{h}\rangle),$$
and this gives \eqref{eq1.8}. $\blacksquare$

Upon integration of \eqref{eq1.8}, one arrives at the {\bf conservation law in the integral form},
$$\int\limits_0^T\langle\sigma u,u\rangle dt+2\Im\langle h_1,h_0\rangle-\left\{\int\limits_0^T\langle\sigma v,v\rangle+2\Im\langle h_1(T),h_0(T)\rangle\right\}$$
$$=2\Im\int\limits_0^T\langle B\dot{h},h\rangle dt.$$
The left-hand side of this equation is the difference between initial ($t=0$) and final ($t=T$) energy of the system, and at the right-hand side there is ``total energy loss'' for the period from 0 to $T$. Use spectral decomposition of the operator $B$,
$$B=\int\limits_{\mathbb{R}}\lambda dE_\lambda,\quad B=B_+-B_-,\quad B_\pm=\pm\int\limits_{\mathbb{R}_\pm}\lambda dE_\lambda\geq0,$$
then the last equation becomes
$$\int\limits_0^T\langle\sigma u,u\rangle dt+2\Im\langle h_1,h_0\rangle+2\Im\int\limits_0^T\langle B_-\dot{h},h\rangle dt$$
$$=\int\limits_0^T\langle\sigma v,v\rangle+2\Im\langle h_1(T),h_0(T)\rangle+2\Im\int\limits_0^T\langle B_t\dot{h},h\rangle dt.$$
The third summand in both sides of the equation describes the losses of the input and output that are caused by dissipation.

Invariancy with regard to coupling is an important property of open systems \cite{23,24}.

\begin{definition}\label{d1.2}
Let
\begin{equation}
\Delta_k=(A_k,H_k,\varphi_k,E_k,\sigma_k)\quad(k=1,2)\label{eq1.9}
\end{equation}
be two local colligations \eqref{eq1.1}, \eqref{eq1.2}, and $\mathcal{F}_{\Delta_k}(B_k)$ be associated with the pairs $\{\Delta_k,B_k\}_1^2$ of the open system \eqref{eq1.6}, \eqref{eq1.7}
\begin{equation}
\mathcal{R}_{\Delta_k}(B_k):\left\{
\begin{array}{lll}
\ddot{h}_k+B_k\dot{h}_k+A_kh_k=\varphi_k^*\sigma_ku_k;\\
\left.h_k\right|_0=h_{k,0},\left.\dot{h}_k\right|_0=h_{k,1}\quad(0\leq t\leq T<\infty,k=1,2);
\end{array}\right.\label{eq1.10}
\end{equation}
$$S_{\Delta_k}:v_k=u_k-i\varphi_kh_k\quad(k=1,2).$$
An open system $\mathcal{F}_{\widetilde{\Delta}}(\widetilde{B})=\{\mathcal{R}_{\widetilde{\Delta}}(\widetilde{B}),S_{\widetilde{\Delta}}\}=\mathcal{F}_{\Delta_2}(B_2)\curlyvee\mathcal{F}_{\Delta_1}(B_1)$ is said to be the {\bf coupling of the open systems} $\{\mathcal{F}_{\Delta_k}(B_k)\}_1^2$ if
\begin{equation}
E=E_1=E_2,\quad\sigma=\sigma_1=\sigma_2,\quad u_2=v_1.\label{eq1.11}
\end{equation}
\end{definition}

Substituting $u_2=v_1-i\varphi_1h_1$ into equation \eqref{eq1.10} for $k=2$, one finds
$$\ddot{h}_2+B_2\dot{h}_2+A_2h_2+i\varphi_2^*\sigma\varphi_1h_1=\varphi_2^*\sigma u_1.$$
Combining this equation with equation \eqref{eq1.10} for $k=1$, one obtains the open system
\begin{equation}
\begin{array}{lll}
{\displaystyle\mathcal{R}_{\widetilde{\Delta}}(\widetilde{B}):\left\{
\begin{array}{lll}
\ddot{\widetilde{h}}+\widetilde{B}\dot{\widetilde{h}}+\widetilde{A}\widetilde{h}=\widetilde{\varphi}^*\sigma u_1;\\
\left.\widetilde{h}\right|_0=\widetilde{h}_0,\left.\dot{\widetilde{h}}\right|_0=\widetilde{h}_1\quad(0\leq t\leq T);
\end{array}\right.}\\
S_{\widetilde{\Delta}}:v_2=u_1-i\widetilde{\varphi}\widetilde{h}
\end{array}\label{eq1.12}
\end{equation}
where $\widetilde{h}=\col[h_1,h_2]$, which corresponds to the traditional \cite{23,24} coupling of colligations $\{\Delta_k\}_1^2$ \eqref{eq1.9}
\begin{equation}
\widetilde{\Delta}=\Delta_2\curlyvee\Delta_1=\left(\widetilde{A}=\left[
\begin{array}{ccc}
A_1&0\\
i\varphi_2^*\sigma\varphi_1&A_2
\end{array}\right],\widetilde{H}=H_1\oplus H_2,\widetilde{\varphi}=\varphi_1+\varphi_2,E,\sigma\right)\label{eq1.13}
\end{equation}
besides,
\begin{equation}
\widetilde{B}=\left[
\begin{array}{ccc}
B_1&0\\
0&B_2
\end{array}\right].\label{eq1.14}
\end{equation}

\begin{remark}\label{r1.2}
Class of open systems $\mathcal{F}_\Delta(B)$ \eqref{eq1.6}, \eqref{eq1.7} is closed with regard to coupling operation. Self-adjointness of the operator $B$ is preserved under coupling \eqref{eq1.14}.
\end{remark}
\vspace{5mm}

{\bf 1.3.} Solution to the Cauchy problem \eqref{eq1.6} is written in terms of the operators $X$ and $Y$ from the factorization of a quadratic operator pencil,
\begin{equation}
L(\lambda)\stackrel{\rm def}{=}\lambda^2I+\lambda B+A=(\lambda I-Y)(\lambda I-X).\label{eq1.15}
\end{equation}
Suppose that decomposition \eqref{eq1.15} exists, where $X$ and $Y$ are linear bounded operators in $H$ and their spectra $\sigma(X)$ and $\sigma(Y)$ do not overlap, $\sigma(X)\cap\sigma(Y)=\emptyset$. Operator $X$ is the right root of the equation
\begin{equation}
X^2+BX+A=0\label{eq1.16}
\end{equation}
since due to \eqref{eq1.15}, $B=-X-Y$, $A=YX$. Denote by $\Gamma$ a closed contour in $\mathbb{C}$ containing $\sigma(Y)$ but not containing $\sigma(X)$ and define the operator $K$ in $H$,
\begin{equation}
K\stackrel{\rm def}{=}\frac1{2\pi i}\int\limits L^{-1}(\zeta)d\zeta.\label{eq1.17}
\end{equation}

{\bf Lemma 1.1.}\label{l1.1}
For the operator $K$ \eqref{eq1.17}, the following equality holds:
\begin{equation}
KY-XK=I.\label{eq1.18}
\end{equation}

P r o o f. Taking into account the identity $Y(\lambda I-Y)^{-1}=\lambda(\lambda I-Y)^{-1}-I$, one obtains
$$KY=\frac1{2\pi i}\int\limits_\Gamma(\zeta I-X)^{-1}(\zeta I-Y)^{-1}Yd\zeta=\frac1{2\pi i}\int\limits_\Gamma\zeta(\zeta I-X)^{-1}(\zeta I-Y)^{-1}d\zeta$$
since $(\lambda I-X)^{-1}$ is holomorphic inside contour $\Gamma$. Again using equation $\zeta(\zeta I-X)^{-1}=X(\zeta I-X)^{-1}+I$, one finds
$$KY=XK+\frac1{2\pi i}\int\limits_\Gamma(\zeta I-Y)^{-1}d\zeta=I,$$
this proves \eqref{eq1.18}.

\begin{corollary}\label{c1.1}
The following relations hold:
\begin{equation}
K(\lambda I-Y)^{-1}-(\lambda I-X)^{-1}K=L^{-1}(\lambda);\quad KY^2+BKY+AK=0\label{eq1.19}
\end{equation}
where $L(\lambda)$ is given by \eqref{eq1.15} and $K$ is given by \eqref{eq1.17}.
\end{corollary}

The first equation in \eqref{eq1.19} follows from \eqref{eq1.18}. To prove the second equation, use $B=-X-Y$, $A=YX$ and $XK=KY-I$,
$$KY^2-YKY-XKY+YXK=KY^2-YKY-(KY-I)Y+Y(KY-I)=0.$$

\begin{remark}\label{r1.3}
Instead of the operator $K$ \eqref{eq1.17}, consider the operator $K_1$,
\begin{equation}
K_1\stackrel{\rm def}{=}\frac1{2\pi i}\int\limits_{L_1}L^{-1}(\zeta)d\zeta\label{eq1.20}
\end{equation}
where $L_1$ is a closed contour in $\mathbb{C}$ containing $\sigma(X)$ and not containing $\sigma(Y)$. Analogously to \eqref{eq1.18}, \eqref{eq1.19}, the following equalities hold:
\begin{equation}
K_1Y-XK_1=-I,\quad K_1(\lambda I-Y)^{-1}-(\lambda I-X)^{-1}K_1=-L^{-1}(\lambda),\label{eq1.21}
\end{equation}
and thus in view of \eqref{eq1.19}
$$(K+K_1)(\lambda I-Y)^{-1}=(\lambda I-X)^{-1}(K+K_1),$$
therefore, taking into account that $\sigma(X)\cap\sigma(Y)=\emptyset$, one obtains that
\begin{equation}
K+K_1=0.\label{eq1.22}
\end{equation}
\end{remark}

The function
\begin{equation}
F_0(t)=e^{tX}g+Ke^{tY}\widetilde{g}\quad(g,\widetilde{g}\in H)\label{eq1.23}
\end{equation}
is a solution to the equation $\ddot{F}+B\dot{F}+AF=0$ which follows from \eqref{eq1.16} and \eqref{eq1.19}.

\begin{remark}\label{r1.4}
Functions
\begin{equation}
G(t,g)=e^{tX}g,\quad\widetilde{G}(t,g)=Ke^{tY}\widetilde{g}\quad(g,\widetilde{g}\in H)\label{eq1.24}
\end{equation}
are linearly independent and form the fundamental system of solutions to the equation $\ddot{y}+B\dot{y}+Ay=0$. So, linear dependency
$$\mu G(t,g)+\nu\widetilde{G}(t,\widetilde{g})=0\quad(\mu,\nu\in\mathbb{C},\mu\not=0,\nu\not=0)$$
after differentiation implies
$$\mu XG(t,g)+\nu X\widetilde{G}(t,g)+\nu e^{tY}\widetilde{g}=0$$
(see \eqref{eq1.18}). Hence it follows $e^{tY}\widetilde{g}=0$, i. e., $\widetilde{g}=0$, this, due to the first equality, gives $g=0$.
\end{remark}

Initial data $F_0(0)=h_0$, $\dot{F}_0(0)=h_1$ for $F_0(t)$ \eqref{eq1.23} in view of \eqref{eq1.18} give the equation system
$$\left\{
\begin{array}{lll}
g+K\widetilde{g}=h_0;\\
Xg+XK\widetilde{g}+\widetilde{g}=h_1;
\end{array}\right.$$
solution of which is
\begin{equation}
\widetilde{g}=h_1-Xh_0,\quad g=h_0-K(h_1-Xh_0).\label{eq1.25}
\end{equation}
So, $F_0(t)$ \eqref{eq1.23} is the solution to Cauchy problem \eqref{eq1.6} for $u=0$ where $g$ and $\widetilde{g}$ are expressed via $h_0$ and $h_1$ by the formulas \eqref{eq1.25}.

Consider the function
$$F_1(t)=\int\limits_0^te^{(t-s)X}f(s)ds+K\int\limits_0^te^{(t-s)Y}\widetilde{f}(s)ds$$
where $f(t)$, $\widetilde{f}(t)$ are vector functions from $H$. Evidently, $F_1(0)=0$. Since
$$\dot{F}_1(t)=X\int\limits_0^te^{(t-s)X}f(s)ds+KY\int\limits_0^te^{(t-s)Y}\widetilde{f}(s)ds+f(t)+K\widetilde{f}(t),$$
then assuming that
\begin{equation}
f(t)+K\widetilde{f}(t)=0\label{eq1.26}
\end{equation}
one obtains that $\dot{F}_1(0)=0$. Differentiating again, one finds that
$$\ddot{F}_1=X^2\int\limits_0^te^{(t-s)X}f(t)dt+KY^2\int\limits_0^te^{(t-s)Y}\widetilde{f}(t)dt+Xf(t)+KY\widetilde{f}(t)$$
and thus $F_1(t)$ is a solution to the equation $\ddot{F}_1(t)+B\dot{F}_1(t)+AF_1(t)=\varphi^*\sigma u(t)$ if only $Xf(t)+KY\widetilde{f}(t)=\varphi^*\sigma u(t)$ (due to \eqref{eq1.16} and \eqref{eq1.19}). Using \eqref{eq1.18} and \eqref{eq1.26}, one obtains that $\widetilde{f}=\varphi^*\sigma u$. Thus, the function
\begin{equation}
F_1(t)=-\int\limits_0^te^{(t-s)x}K\varphi^*\sigma u(s)ds+K\int\limits_0^te^{(t-s)Y}\varphi^*\sigma u(s)ds\label{eq1.27}
\end{equation}
is the solution to the Cauchy problem \eqref{eq1.6} for zero initial data, $h_0=h_1=0$.

\begin{theorem}\label{t1.2}
If factorization \eqref{eq1.15} where $X$ and $Y$ are linear operators in $H$ and $\sigma(X)\cap\sigma(Y)=\emptyset$, then the solution to Cauchy problem \eqref{eq1.6} is $h=F_0+F_1$ where $F_0$ and $F_1$ are given by \eqref{eq1.23} and \eqref{eq1.27} respectively, the operator $K$ is given by formula \eqref{eq1.17} and $g$ and $\widetilde{g}$ are expressed via $h_0$ and $h_1$ by equations \eqref{eq1.25}.
\end{theorem}
\vspace{5mm}

{\bf 1.4.} To the input $u$ of open system \eqref{eq1.6}, \eqref{eq1.7}, give a plain wave signal $u=e^{\lambda t}u_0$ ($\lambda\in\mathbb{C}$, $u_0\in E$ do not depend on $t$) and let $\lambda\not\in\sigma(X)\cup\sigma(Y)$. Since
$$\int\limits_0^te^{(t-s)X}e^{\lambda s}ds=(\lambda I-X)^{-1}\left(e^{\lambda t}I-e^{tX}\right),$$
the function $F_1(\lambda)$ \eqref{eq1.27} in this case is
$$F_1(t)=e^{tX}(\lambda I-X)^{-1}\varphi^*\sigma u_0+e^{\lambda t}\{K(\lambda I-Y)^{-1}-(\lambda I-X)^{-1}K\}\varphi^*\sigma u_0$$
$$-Ke^{tY}(\lambda I-Y)^{-1}\varphi^*\sigma u_0$$
or taking \eqref{eq1.19} into account, one has that
$$F_1(t)=e^{\lambda t}L^{-1}(\lambda)\varphi^*\sigma u_0+e^{tX}(\lambda I-X)^{-1}\varphi^*\sigma u_0-Ke^{tY}(\lambda I-Y)^{-1}\varphi^*\sigma u_0.$$
So, solution $h=h(t)$ of the Cauchy problem \eqref{eq1.6} for $u=e^{\lambda t}u_0$ is
$$h(t)=e^{tX}\{h_0-K(h_1-Xh_0)+(\lambda I-X)^{-1}\varphi^*\sigma u_0\}+e^{\lambda t}L^{-1}(\lambda)\varphi^*\sigma u_0$$
\begin{equation}
+Ke^{tY}\{h_1-Xh_0-(\lambda I-Y)^{-1}\varphi^*\sigma u_0\}.\label{eq1.28}
\end{equation}
Assuming that internal state $h(t)$ of the open system \eqref{eq1.6}, \eqref{eq1.7} ``oscillates with the same frequency $\lambda$'' as an input signal $u=e^{\lambda t}u_0$, set
\begin{equation}
h(t)=e^{\lambda t}L^{-1}(\lambda)\varphi^*\sigma u_0\label{eq1.29}
\end{equation}
and
$$h_0-K(h_1-Xh_0)+(\lambda I-X)^{-1}\varphi^*\sigma u_0=0;\quad h_1-Xh_0-(\lambda I-Y)^{-1}\varphi^*\sigma u_0=0,$$
this, due to \eqref{eq1.19}, gives
\begin{equation}
h_0=L^{-1}(\lambda)\varphi^*\sigma u_0,\quad h_1=\lambda L^{-1}(\lambda)\varphi^*\sigma u_0.\label{eq1.30}
\end{equation}
Equations \eqref{eq1.30} completely correlate with formula \eqref{eq1.29}. Assuming that the output signal $v$ of the open system $\mathcal{F}_\Delta(B)$ \eqref{eq1.6}, \eqref{eq1.7} is $v=e^{\lambda t}v_0$, from \eqref{eq1.7}, due to \eqref{eq1.29}, one finds that
\begin{equation}
v_0=S_\Delta(\lambda,B)u_0\label{eq1.31}
\end{equation}
where $S_\Delta(\lambda,B)$ is
\begin{equation}
S_\Delta(\lambda,B)=I-i\varphi L^{-1}(\lambda)\varphi^*\sigma\label{eq1.32}
\end{equation}
and is said to be the {\bf characteristic function of the pair} $\{\Delta,B\}$ or {\bf characteristic function of the pencil} $L(\lambda)$ \eqref{eq1.15}.

The equality
$$\sigma-S_\Delta^*(w,B)\sigma S_\Delta(\lambda,B)=i\sigma\varphi\{L^{-1}(\lambda)-\left(L^*(w)\right)^{-1}+\left(L^*(w)\right)^{-1}i\varphi^*\sigma\varphi L^{-1}(\lambda)\}\varphi^*\sigma$$
$$=i\sigma\varphi\left(L^*(w)\right)^{-1}\left[\left(\overline{w}^2-\lambda^2\right)I+(\overline{w}-\lambda)B+A^*-A+i\varphi^*\sigma\varphi\right]L^{-1}(\lambda)\varphi^*\sigma,$$
due to $B=B^*$ and \eqref{eq1.2}, implies
\begin{equation}
\frac i{\lambda-\overline{w}}\{\sigma-S_\Delta^*(w,B)\sigma S_\Delta(\lambda,B)\}=\sigma\varphi\left(L^*(w)\right)^{-1}[(\lambda+\overline{w})I+B]L^{-1}(\lambda)\varphi^*\sigma.\label{eq1.33}
\end{equation}
This {\bf metric relation} for characteristic function $S_\Delta(\lambda,B)$ \eqref{eq1.32} follows also from the conservation law if one takes into account \eqref{eq1.29} and \eqref{eq1.31}.

The coupling rules \eqref{eq1.13}, \eqref{eq1.14} imply that characteristic functions
\begin{equation}
S_{\Delta_k}(\lambda,B_k)=I-i\varphi_kL_k^{-1}(\lambda)\varphi_k^*\sigma\quad(k=1,2)\label{eq1.34}
\end{equation}
of pairs $\{\Delta_k,B_k\}$ are multiplied under coupling,
\begin{equation}
S_{\widetilde{\Delta}}(\lambda,\widetilde{B})=S_{\Delta_2}(\lambda,B_2)S_{\Delta_1}(\lambda,B_1)\label{eq1.35}
\end{equation}
(here $\widetilde{\Delta}$ and $\widetilde{B}$ are given by \eqref{eq1.13} and \eqref{eq1.14}).

Let the pencils $L_k(\lambda)$ have factorizations
\begin{equation}
L_k(\lambda)=\lambda^2I+\lambda B_k+A_k=(\lambda I-Y_k)(\lambda I-X_k)\quad(k=1,2)\label{eq1.36}
\end{equation}
where $\sigma(X_k)\cap\sigma(Y_k)=\emptyset$ ($k=1$, 2) and
\begin{equation}
-B_k=X_k+Y_k,\quad A_k=Y_kX_k\quad(k=1,2).\label{eq1.37}
\end{equation}
Give method for constructing ``roots'' of the pencil
\begin{equation}
\widetilde{L}(\lambda)=\lambda^2I+\lambda\widetilde{B}+\widetilde{A}=(\lambda I-\widetilde{Y})(\lambda I-\widetilde{X})\label{eq1.38}
\end{equation}
($\widetilde{A}$, $\widetilde{B}$ are given by \eqref{eq1.13}, \eqref{eq1.14}) by the roots $\{X_k,Y_k\}_1^2$ of the pencils $\{L_k(\lambda)\}_1^2$ \eqref{eq1.36}. Let
\begin{equation}
\widetilde{X}=\left[
\begin{array}{ccc}
X_1&0\\
\gamma&X_2
\end{array}\right],\quad\widetilde{Y}=\left[
\begin{array}{ccc}
Y_1&0\\
-\gamma&Y_2
\end{array}\right]\quad(\gamma:H_1\rightarrow H_2),\label{eq1.39}
\end{equation}
then the condition $\widetilde{X}+\widetilde{Y}=-\widetilde{B}$ holds due to \eqref{eq1.14}, \eqref{eq1.37}. Since
$$\widetilde{Y}\widetilde{X}=\left[
\begin{array}{ccc}
Y_1X_1&0\\
Y_2\gamma-\gamma X_1&Y_2X_2
\end{array}\right],$$
equations \eqref{eq1.13} and \eqref{eq1.37} imply that $\widetilde{A}=\widetilde{Y}\widetilde{X}$ if only
\begin{equation}
Y_2\gamma-\gamma X_1=i\varphi_2^*\sigma\varphi_1.\label{eq1.40}
\end{equation}
Equation \eqref{eq1.40} is solvable \cite{22} under condition that
\begin{equation}
\sigma(Y_2)\cap\sigma(X_1)=\emptyset.\label{eq1.41}
\end{equation}
So, if $\gamma$ is a solution to equation \eqref{eq1.40}, then $\widetilde{X}$, $\widetilde{Y}$ \eqref{eq1.39} are roots of the pencil $\widetilde{L}(\lambda)$ \eqref{eq1.38}. Probably, this method of construction of roots $\widetilde{X}$, $\widetilde{X}$, $\widetilde{Y}$ by $\{X_k,Y_k\}$ is not unique.
\vspace{5mm}

{\bf 1.5.} Substitute $u=v+i\varphi h$ (see \eqref{eq1.7}) into equation \eqref{eq1.6}, then, taking into account \eqref{eq1.2}, one has
$$\ddot{h}+B\dot{h}+A^*h=\varphi^*\sigma v.$$
Half-sum of this equation and \eqref{eq1.6} gives \cite{23,24} {\bf diagonal of the open system} $\mathcal{F}_d(B)=\{\mathcal{R}_d(B),S_d\}$ where
\begin{equation}
\mathcal{R}_d(B):\left\{
\begin{array}{lll}
\ddot{h}+B\dot{h}+A_Rh=\varphi^*\sigma u_d;\\
\left.h\right|_0=h_0,\left.\dot{h}\right|_0=h_1\,(0<t<T);
\end{array}\right.\label{eq1.42}
\end{equation}
$$S_d:v_d=\frac i2\varphi h,$$
here
\begin{equation}
2A_R=A+A^*,\quad 2u_d=u+v,\quad 2v_d=u-v.\label{eq1.43}
\end{equation}
Define the operator-valued function
\begin{equation}
V(\lambda,B)\stackrel{\rm def}{=}\varphi L_R^{-1}(\lambda)\varphi^*\quad(L_R(\lambda)=\lambda^2I+\lambda B+A_R),\label{eq1.44}
\end{equation}
then \eqref{eq1.43}, due to \eqref{eq1.31}, yields that
$$I-S_\Delta(\lambda,B)=\frac i2V(\lambda,B)\sigma(I+S_\Delta(\lambda,B)),$$
i. e., $V(\lambda,B)$ \eqref{eq1.44} is expressed via $S_\Delta(\lambda,B)$ by the fractional linear transformation,
\begin{equation}
V(\lambda,B)=2i(S_\Delta(\lambda,B)-I)(S_\Delta(\lambda,B)+I)^{-1}\sigma^{-1}.\label{eq1.45}
\end{equation}
Notice that the function $V(\lambda,B)$ corresponds to the self-adjoint pencil $L_R^*(\lambda)=L_R(\overline{\lambda})$ and
\begin{equation}
\frac{V(\lambda,B)-V^*(w,B)}{\lambda-\overline{w}}=-\varphi(L_R^*(w))^{-1}[(\lambda+\overline{w})I+B]L_R^{-1}(\lambda)\varphi^*.\label{eq1.46}
\end{equation}

\section{Functional model of a self-adjoint operator}\label{s2}

{\bf 2.1.} Let $T$ be a self-adjoint operator acting in a separable Hilbert space $H$ and $E_t$ ($t\in\mathbb{R}$) be its decomposition of identity \cite{25} -- \cite{26}. Consider the subspace $L_1$ in $H$,
\begin{equation}
L_1\stackrel{\rm def}{=}\span\{E_\Delta f_1:\Delta\in\mathbb{R}\},\label{eq2.1}
\end{equation}
here $\Delta$ runs through all the intervals from $\mathbb{R}$ and $f_1$ is a fixed vector from a countable (separable) dense set in $H$. For $L_1$ \eqref{eq2.1} and $H_1=H\ominus L_1$, $E_tL_1\subseteq L_1$, $E_tH_1\subseteq H_1$ for all $t\in\mathbb{R}$. Subset $H_1$ is separable. Analogously to \eqref{eq2.1}, define the subspace $L_2$ in $H_1$,
$$L_2\stackrel{\rm def}{=}\span\{E_\Delta f_2:\Delta\in\mathbb{R}\}$$
where, as before, $\Delta$ runs through the set of all intervals from $\mathbb{R}$ and $f_2$ is a fixed vector from a dense set in $H_1$. Countable repetition of this procedure (due to separability of $H$) gives
\begin{equation}
H=\sum\limits_1^\infty\oplus L_k=\span\{E_\Delta g:g\in G,\Delta\in\mathbb{R}\}\label{eq2.2}
\end{equation}
where $G$ is a {\bf generating subspace} \cite{25} of the operator $T$,
\begin{equation}
G\stackrel{\rm def}{=}\span\{f_k:k\in\mathbb{N}\}.\label{eq2.3}
\end{equation}

\begin{remark}\label{r2.1} The fact that orthogonal sum \eqref{eq2.2} of the subspaces of \eqref{eq2.1} type exhausts the whole $H$ follows from separability of $H$. Really, let ${\displaystyle H_1=\sum\limits_1^\infty\oplus L_k}$ where $L_k$ are constructed in above way and $H_0=H\ominus H_1$. Then $E_tH_s\subset H_s$ ($s=0$, $1$) and $H_0$ is separable. Choosing $f_0\in H_0$ from the countable dense set in $H_0$, construct $L_0=\span\{E_\Delta f_0:\Delta\in\mathbb{R}\}$. $L_0\subset H_1$ since $H_1$ contains all the subspaces of such type, therefore $L_0=\{0\}$.
\end{remark}

In the case of boundedness of $T$,
\begin{equation}
H=\span\{T^ng:g\in G;n\in\mathbb{Z}_+\}.\label{eq2.4}
\end{equation}
Denote by $\mathfrak{N}$ a Hilbert space ($\dim\mathfrak{N}=\dim G$) and by $\psi$ denote a linear bounded operator from $H$ onto $\mathfrak{N}$. In the case of $\mathfrak{N}=G$, one can put $\psi=P_G$ where $P_G$ is an orthogonal projection onto $G$. Define in $\mathfrak{N}$ a non-decreasing non-negative operator-valued function $F(t)\stackrel{\rm def}{=}\psi E_t\psi^*$ ($t\in\mathbb{R}$, $E_t$ is a resolution of the identity of $T$). Define the Hilbert space of vector functions
\begin{equation}
L_{\mathbb{R}}^2(\mathfrak{N},dF)\stackrel{\rm def}{=}\left\{f(x)\in\mathfrak{N}:\int\limits_{\mathbb{R}}\langle dF(x)f(x),f(x)\rangle_{\mathfrak{N}}<\infty\right\}\label{eq2.5}
\end{equation}
obtained as a result of closure of a class of continuous $\mathfrak{N}$-valued functions with consequent factorization by the metric kernel \cite{19,26}. Denote by $Q$  the self-adjoint operator of multiplication by independent variable in $L_{\mathbb{R}}^2(\mathfrak{N},dF)$,
\begin{equation}
(Qf)(x)\stackrel{\rm def}{=}xf(x)\quad(f\in L_{\mathbb{R}}^2(\mathfrak{N},dF)),\label{eq2.6}
\end{equation}
domain of which is
\begin{equation}
\mathfrak{D}_Q\stackrel{\rm def}{=}\{f\in L_{\mathbb{R}}^2(\mathfrak{N},dF):xf(x)\in L_{\mathbb{R}}^2(\mathfrak{N},dF)\}.\label{eq2.7}
\end{equation}
Consider the linear operator $V:$ $L_{\mathbb{R}}^2(\mathfrak{N},dF)\rightarrow H$,
\begin{equation}
Vf(.)=f;\quad f\stackrel{\rm def}{=}\int\limits_{\mathbb{R}}dE_x\psi^*f(x).\label{eq2.8}
\end{equation}

\begin{theorem}\label{t2.1}
Let $G$ \eqref{eq2.3} be the generating subspace \eqref{eq2.2} of a self-adjoint operator $T$ acting in a separable Hilbert space $H$ and $\psi:$ $H\rightarrow\mathfrak{N}$ be a linear bounded operator such that $\Ker\psi=H\oplus G$ and $\psi^*$ be invertible on $\mathfrak{N}$. Then the operator $V$ is a unitary isomorphism between $L_{\mathbb{R}}^2(\mathfrak{N},dF)$ \eqref{eq2.5} ($F(x)=\psi E_x\psi^*$, $E_x$ is the resolution of the identity of $T$) and $H$, besides,
\begin{equation}
VQ=TV\label{eq2.9}
\end{equation}
where $Q$ is given by \eqref{eq2.6}, \eqref{eq2.7}.
\end{theorem}

P r o o f. Equation \eqref{eq2.8} implies that
$$\psi E_xf=\int\limits_{-\infty}^xdF(t)f(t),$$
therefore
$$\langle f,f\rangle=\left\langle f,\int\limits_{\mathbb{R}}dE_x\psi^*f(x)\right\rangle=\int\limits_{\mathbb{R}}\langle d(\psi E_xf),f(x)\rangle=\int\limits_{\mathbb{R}}\langle dF(x)f(x),f(x)\rangle,$$
this proves that $V$ \eqref{eq2.8} is isometric. Image of $V$ is dense in $H$ since vectors of the form $E_\Delta\psi^*g$ ($\Delta\in\mathbb{R}$, $g\in\mathfrak{N}$) belong to it and thus $V$ is unitary. Relation \eqref{eq2.9} follows from the equation
$$TVf(.)=T\int\limits_{\mathbb{R}}dE_x\psi^*f(x)=\int\limits_{\mathbb{R}}dE_x\psi^*xf(x)=VQf(.),$$
besides,
$$\|Tf\|_H^2=\int\limits_{\mathbb{R}}x^2\langle dF(x)f(x),f(x)\rangle=\|Qf(x)\|_{L^2}^2,$$
for all $f$ of \eqref{eq2.8} type. Thus $V$ maps domain $\mathfrak{D}_T$ of the operator $T$ onto domain $\mathfrak{D}_Q$ \eqref{eq2.7} of the operator $Q$. $\blacksquare$

Description of the commutant of the operator $Q$ \eqref{eq2.6}, \eqref{eq2.7} is given by the theorem \cite{26}.

\begin{theorem}\label{t2.2}
Any bounded operator $N$ in $L_{\mathbb{R}}^2(\mathfrak{N},dF)$ commutating with $Q$ \eqref{eq2.6}, \eqref{eq2.7} is
\begin{equation}
(Nf)(x)=n(x)f(x)\label{eq2.10}
\end{equation}
where $n(x)$ is an operator-valued function in $\mathfrak{N}$, besides,
\begin{equation}
\|N\|=\sup\limits_{x\in\mathbb{R}}(F)\|n(x)\|_{\mathfrak{N}}.\label{eq2.11}
\end{equation}
\end{theorem}

Commutativity of $N$ and $Q$ is defined as commutativity of $N$ and $R_Q(\lambda)=(Q-\lambda I)^{-1}$ at points of regularity $\lambda\in\mathbb{C}$. Symbol $\sup(F)$ means \cite{26} that $\sup$ is taken almost everywhere by the measure $\langle dFf,f\rangle$ for all $f\in\mathfrak{N}$.

P r o o f. $[N,Q]=0$ implies that the operator $\widetilde{N}=VNV^*$ ($V$ is given by \eqref{eq2.8}) commutes with $T$, therefore $\widetilde{N}E_x=E_x\widetilde{N}$ ($\forall x$, \cite{25}), thus
\begin{equation}
\widetilde{N}f=\int\limits_{\mathbb{R}}dE_x\widetilde{N}\psi^*f(x).\label{eq2.12}
\end{equation}
Denote by $\{e_k\}_1^\infty$ an orthonormal basis in $\mathfrak{N}$. To the vector $\widetilde{N}\psi^*e_k\in H$ there corresponds the function $h_k(.)\in L_{\mathbb{R}}^2(\mathfrak{N},dF)$,
$$\widetilde{N}\psi^*e_k=\int\limits_{\mathbb{R}}dE_x\psi^*h_k(x)\quad(k\in\mathbb{N}),$$
due to unitarity of $V$ \eqref{eq2.8}. ${\displaystyle I=\sum\limits_k\langle.,e_k\rangle e_k}$ yields
\begin{equation}
\widetilde{N}\psi^*=\int\limits_\mathbb{R}dF_x\psi^*n(x),\label{eq2.13}
\end{equation}
where $n(x)$ is an operator-valued function in $\mathfrak{N}$,
$$n(x)\stackrel{\rm def}{=}\sum\limits_k\langle.,e_k\rangle h_k(x).$$
Using \eqref{eq2.12}, \eqref{eq2.13}, one obtains
$$\widetilde{N}f=\int\limits_{\mathbb{R}}dF(x)\psi^*n(x)f(x),$$
hence \eqref{eq2.10} follows. Equation \eqref{eq2.11} is proved in a standard way \cite{25}. $\blacksquare$
\vspace{5mm}

{\bf 2.2} Let $K$ be a finite-dimensional self-adjoint operator in $H$ and $r=\rank K<\infty$. Rewrite $K$ as $K=|K|^{1/2}\sign K\cdot|K|^{1/2}$ ($|K|^{1/2}$, $\sign K$ are understood in the sense of spectral decompositions \cite{25,26}), then
\begin{equation}
K=\sum\limits_{\alpha,\beta=1}^r\langle.,g_\alpha\rangle\sigma_{\alpha,\beta}g_\beta,\label{eq2.14}
\end{equation}
here $g_\alpha=|K|^{1/2}e_\alpha$; $\{e_\alpha\}_1^2$ is an orthonormal basis in $KH$; $\sigma_{\alpha,\beta}$ are matrix elements of $\sign K$ in the basis $\{e_\alpha\}_1^r$. Equation \eqref{eq2.8} implies that
$$g_\alpha=\int\limits_\mathbb{R}dE_x\psi^*g_\alpha(x)\quad(1\leq\alpha\leq r)$$
where $g_\alpha(.)\in L_{\mathbb{R}}^2(\mathfrak{N},dF)$ and $E_x$ is the resolution of the identity of $T$. For $K$ \eqref{eq2.14}, the following representation holds,
$$Kf=\int\limits_\mathbb{R}dE_x\psi^*\sum\limits_{\alpha,\beta}\int\limits_\mathbb{R}\langle dF(t)f(t),g_\alpha(t)\rangle\sigma_{\alpha,\beta}g_\beta(x)=\int\limits_\mathbb{R}dE_x\psi^*\int\limits_\mathbb{R}K(x,t)dF(t)f(t)$$
where $K(x,t)$ is an operator-valued function in $\mathfrak{N}$,
$$K(x,t)\stackrel{\rm def}{=}\sum\limits_{\alpha,\beta=1}^r\langle.,g_\alpha(t)\rangle\sigma_{\alpha,\beta}g_\beta(x).$$
So, the operator $\widetilde{K}=V^*KV$ in the space $L_{\mathbb{R}}^2(\mathfrak{N},dF)$ is
\begin{equation}
(\widetilde{K}f)(x)=\int\limits_\mathbb{R}K(x,t)dF(t)f(t)\quad(f(.)\in L_{\mathbb{R}}^2(\mathfrak{N},dF)),\label{eq2.15}
\end{equation}
besides, $K^*(x,t)=K(t,x)$ in view of self-adjointness of $K$.

\begin{theorem}\label{t2.3}
Let $K$ be a bounded self-adjoint operator in $H$, then the operator $\widetilde{K}=V^*KV$ ($V$ is given by \eqref{eq2.8}) in the space $L_{\mathbb{R}}^2(\mathfrak{N},dF)$ is an integral operator \eqref{eq2.15} where kernel $K(x,t)$ is an operator-valued function in $\mathfrak{N}$; $K^*(x,t)=K(t,x)$, besides, $K(x,t)$ is such that $\widetilde{K}$ \eqref{eq2.15} is bounded in $L_{\mathbb{R}}^2(\mathfrak{N},dF)$.
\end{theorem}

So, an arbitrary bounded self-adjoint operator $K$ in the spectral representation of the operator $\widetilde{T}=Q$ \eqref{eq2.6}, \eqref{eq2.7} is realized by  the integral operator $\widetilde{K}$ \eqref{eq2.15} in the space $L_{\mathbb{R}}^2(\mathfrak{N},dF)$ \eqref{eq2.5}.

\section{Model representations of quadratic pencils}

{\bf 3.1.} Consider the operator-valued quadratic pencil \eqref{eq1.15},
\begin{equation}
L(\lambda)=\lambda^2I+\lambda B+A\quad(\lambda\in\mathbb{C})\label{eq3.1}
\end{equation}
where $B$ and $A$ are bounded linear operators in a Hilbert space $H$, besides, $B=B^*$, and for $A$ \eqref{eq1.2} holds. Polynomial operator pencils are studied by different methods \cite{1,2}, \cite{6} -- \cite{16}, \cite{29}, \cite{31} -- \cite{34}, one of which consists in the use of factorization
\begin{equation}
L(\lambda)=(\lambda I-Y)(\lambda I-X)\label{eq3.2}
\end{equation}
where $X$, $Y$ are linear bounded operators in $H$ such that
\begin{equation}
\left\{
\begin{array}{ccc}
X+Y=-B;\\
YX=A.
\end{array}\right.\label{eq3.3}
\end{equation}
Hence it follows that $X(Y)$ is a right (left) operator root of the equation
\begin{equation}
X^2+BX+A=0\quad(Y^2+YB+A=0).\label{eq3.4}
\end{equation}

\begin{definition}\label{d3.1}
Operator $X$ from factorization \eqref{eq3.2} is said to be the {\bf spectral root} \cite{1} of pencil $L(\lambda)$ \eqref{eq3.2} if $\sigma(X)\cap\sigma(Y)=\emptyset$ where $\sigma(X)$ and $\sigma(Y)$ are spectra of the operators $X$ and $Y$. Decomposition \eqref{eq3.2} in this case is said to be the {\bf spectral factorization}.
\end{definition}

The following statement gives sufficient conditions for existence of a spectral root \cite{1}.

\begin{lemma}\label{l3.1}
If operator $B$ of the quadratic pencil \eqref{eq3.1} is invertible and $4\|B^{-1}\|\cdot\|B^{-1}A\|<1$, then the spectral root $X$ (spectral factorization \eqref{eq3.2}) exists.
\end{lemma}

Since $L(\lambda+\varepsilon)=\lambda^2I+\lambda B_\varepsilon+A_\varepsilon$ where $B_\varepsilon=B+2\varepsilon I$, $A_\varepsilon=A+\varepsilon B+\varepsilon^2I$, then one always can choose such $\varepsilon\in\mathbb{C}$ that $B_\varepsilon$ is invertible.
\vspace{5mm}

{\bf 3.2.} Supposing that spectral factorization \eqref{eq3.2} holds, rewrite root $X$ as
\begin{equation}
X=C+iD\quad(2C=X+X^*,\,2iD=X-X^*),\label{eq3.5}
\end{equation}
then
$$Y=-B-C-iD,$$
and thus the operator $A$ \eqref{eq3.3} is
$$A=-(B+C+iD)(C+iD)=-BC-C^2+D^2+i[DC+(B+C)D].$$
Equation \eqref{eq1.2} implies
\begin{equation}
[C,B]-i\{D,B+2C\}=i\varphi^*\sigma\varphi^*\sigma\varphi\label{eq3.6}
\end{equation}
where $[X,Y]=XY-YX$ is commutator and $[X,Y]=XY+YX$ is anti-commutator of the operators $X$ and $Y$.

\begin{remark}\label{r3.1}
${\rm(i)}$ If $D=0$ (solution $X$, $Y$ to system \eqref{eq3.3} is sought in the class of self-adjoint operators), then \eqref{eq3.6} implies
\begin{equation}
[X,B]=i\varphi^*\sigma\varphi\quad(X=C).\label{eq3.7}
\end{equation}
This case is said to be {\bf hyperbolic}.

${\rm(ii)}$ When $C=0$ ({\bf conditionally elliptic} case), then $X=iD$ \eqref{eq3.5} is skew Hermitian ($X^*=-X$) and $Y=-B-iD$ and
\begin{equation}
\{D,B\}=-\varphi^*\sigma\varphi.\label{eq3.8}
\end{equation}

${\rm(iii)}$ Pure {\bf elliptic} case $X=iD$, $Y=-iD$ gives $B=0$, this case is of little interest.
\end{remark}

In terms of the operator
\begin{equation}
K\stackrel{\rm def}{=}B+2C,\label{eq3.9}
\end{equation}
equation \eqref{eq3.6} becomes
\begin{equation}
\frac12[K,B]-i\{D,K\}=i\varphi^*\sigma\varphi.\label{eq3.10}
\end{equation}
Finding the pairs of self-adjoint operators $K$, $D$ satisfying \eqref{eq3.10} (or solutions $C$, $D$ to equation \eqref{eq3.6}) consists in the following step-by-step procedure.

(i) From the equation
\begin{equation}
\frac12[K,B]=i(\varphi^*\sigma\varphi+T)\label{eq3.11}
\end{equation}
one finds the operator $K$ supposing that $B$ $(=B^*)$, $\varphi^*\sigma\varphi$ and $T=T^*$ are given (here $T$ is an arbitrary bounded self-adjoint operator).

(ii) Knowing $K$ and $T$ (from (i)), define $D$ as solution to the equation
\begin{equation}
\{D,K\}=T.\label{eq3.12}
\end{equation}

Two extreme cases are possible:
\begin{equation}
\begin{array}{lll}
{\rm(a)}\,[K,B]=2i\varphi^*\sigma\varphi;\,\{D,K\}=0\,(T=0);\\
{\rm(b)}\,[K,B]=0;\,\{D,K\}=-\varphi^*\sigma\varphi\,(T=-\varphi^*\sigma\varphi).
\end{array}\label{eq3.13}
\end{equation}

\begin{remark}\label{r3.2}
Self-adjoint operator $T$ is a free parameter of the step-by-step procedure. Equations \eqref{eq3.7}, \eqref{eq3.8} coincide with the equations of step-by-step procedure \eqref{eq3.11}, \eqref{eq3.12} and with relations ${\rm (a)}$, $\rm{(b)}$ \eqref{eq3.13}. Methods for solving the operator equations \eqref{eq3.7}, \eqref{eq3.8} are given below.
\end{remark}
\vspace{5mm}

{\bf 3.3.} Consider the hyperbolic case when $D=0$ and $X$, $Y$ are self-adjoint operators, \eqref{eq3.7} takes place and $Y=-B-X$.

A self-adjoint operator acting in a separable Hilbert space $H$ is realized (Theorem \ref{t2.1}) by the operator \eqref{eq2.6} of multiplication by independent variable,
\begin{equation}
(\widetilde{T}f)(x)=xf(x)\quad(f\in L_{\mathbb{R}}^2(\mathfrak{N},dF)\label{eq3.14}
\end{equation}
where $\widetilde{T}=V^*TV$ and $V$ is given by \eqref{eq2.8}. An arbitrary self-adjoint bounded operator $K$ ($h\rightarrow H$) in the spectral representation of the operator $\widetilde{T}$ \eqref{eq3.14} is given (Theorem \ref{t2.3}) by the integral operator \eqref{eq2.15}
\begin{equation}
(\widetilde{K}f)=\int\limits_\mathbb{R}K(x,t)dF(t)f(t)\quad(f\in L_{\mathbb{R}}^2(\mathfrak{N},dF),\label{eq3.15}
\end{equation}
besides, $\widetilde{K}=V^*KV$; $K(x,t)$ is an operator-valued function in $\mathfrak{N}$ and $K^*(x,t)=K(t,x)$.

\begin{definition}\label{d3.2}
A self-adjoint operator $K$ in $H$ is said to be $T$-{\bf bounded} if its kernel $K(x,t)$ from \eqref{eq3.15} in the spectral representation of $\widetilde{T}$ \eqref{eq3.14} satisfies the condition
\begin{equation}
\sup\limits_{x,t}(F)|\langle K(x,t)f,g\rangle|<C\|f\|\cdot\|g\|\quad(\forall f,g\in\mathfrak{N})\label{eq3.16}
\end{equation}
and the number $C$ does not depend on $f$ and $g$.
\end{definition}

\begin{lemma}\label{l3.2}
Let $K$ and $T$ be bounded self-adjoint operators in a separable Hilbert space $H$ and $K$ be $T$-bounded, then the equation $XT-TX=iK$ always has bounded self-adjoint solution $X$.

In the spectral representation \eqref{eq3.14} of the operator $\widetilde{T}$ solution $\widetilde{X}$ of the equation
\begin{equation}
\widetilde{X}\widetilde{T}-\widetilde{T}\widetilde{X}=i\widetilde{K},\label{eq3.17}
\end{equation}
here $\widetilde{K}$ is given by \eqref{eq3.15}, is given by the formula
\begin{equation}
(\widetilde{X}f)(x)=n(x)f(x)+i\int\limits_\mathbb{R}\frac{K(x,t)}{t-x}dF(t)f(t)\quad(f\in L_\mathbb{R}^2(\mathfrak{N},dF)\label{eq3.18}
\end{equation}
here $n(x)$, $K(x,t)$ are operator-valued functions in $\mathfrak{N}$; $K^*(x,t)=K(t,x)$ and \eqref{eq3.16} takes place; $dFn=n^*dF$ and $\sup\limits_x|\langle n(x)f,g\rangle|<C\|f\|\cdot\|g\|$ ($\forall f$, $g\in\mathfrak{N}$).
\end{lemma}

P r o o f. The operator
$$(\widetilde{X}_0f)(x)\stackrel{\rm def}{=}i\int\limits_{\mathbb{R}}\frac{K(x,t)}{t-x}dF(t)f(t)\quad(f\in L_\mathbb{R}^2(\mathfrak{N},dF)$$
is self-adjoint and bounded (due to \eqref{eq3.16}), besides, $\widetilde{X}_0\widetilde{T}-\widetilde{T}\widetilde{X}_0=i\widetilde{K}$. If $\widetilde{X}$ is the self-adjoint solution of equation \eqref{eq3.17}, then the operator $\widetilde{X}-\widetilde{X}_0$ commutes with $\widetilde{T}$, this, due to Theorem \ref{t2.2}, gives representation \eqref{eq3.18}. $\blacksquare$

\begin{theorem}\label{t3.1}
Let linear bounded operators $A$ and $B$ be defined in a separable Hilbert space $H$ such that

{\rm (i)} $B=B^*$ and $B$ is invertible;

{\rm (ii)} equation \eqref{eq1.2} holds for $A$ and $\varphi^*\sigma\varphi$ is $B$-bounded;

{\rm (iii)} $4\|B^{-1}\|\cdot\|B^{-1}A\|<1$.

Then the quadratic pencil $L(\lambda)$ \eqref{eq3.1} has a spectral root. In the spectral representation of the operator $\widetilde{B}$ (Theorem \ref{t2.1}) the operator $\widetilde{X}$ is
\begin{equation}
(\widetilde{X}f)(x)=n(x)f(x)+i\int\limits_\mathbb{R}\frac{K(x,t)}{t-x}dF(t)f(t)\quad(f\in L_\mathbb{R}^2(\mathfrak{N},dF)),\label{eq3.19}
\end{equation}
here $n(x)$ and $K(x,t)$ have the same properties as in Lemma \ref{l3.2}, kernel $K(x,t)$ corresponds to the integral representation \eqref{eq3.15} of the operator $\widetilde{K}=\widetilde{\varphi^*\sigma\varphi}$.
\end{theorem}

Formula \eqref{eq3.19} follows from \eqref{eq3.18} where $T=B$ and $K=\varphi^*\sigma\varphi$.

Operator $\widetilde{X}$ \eqref{eq3.19} is the sum of the operator of multiplication by the operator-valued function $n(x)$ and the integral operator which is a generalization of a well-known {\bf Hilbert transform}. Functional realization of the root $Y$ is
\begin{equation}
(\widetilde{Y}f)(x)=-(x+n(x))f(x)-i\int\limits_\mathbb{R}\frac{K(x,t)}{t-x}dF(t)f(t)\quad(f\in L_\mathbb{R}^2(\mathfrak{N},dF))\label{eq3.20}
\end{equation}

So, in terms of spectral decomposition of the operator $B$ ($\widetilde{B}=\widetilde{T}$ \eqref{eq3.14}) the roots $\widetilde{X}$ and $\widetilde{Y}$ of the quadratic pencil $L(\lambda)$ \eqref{eq3.2} are given by \eqref{eq3.19} and \eqref{eq3.20} respectively.

Sometimes it is more convenient to write realizations corresponding to the spectral decomposition of $X$,
\begin{equation}
(\widetilde{X}f)(x)=xf(x)\quad(f\in L_\mathbb{R}^2(\mathfrak{N},dF)).\label{eq3.21}
\end{equation}
Then the operator $B$ in representation \eqref{eq3.21}, due to Lemma \ref{l3.2}, is
\begin{equation}
(\widetilde{B}f)(x)=b(x)f(x)-i\int\limits_\mathbb{R}\frac{K(x,t)}{t-x}dF(t)f(t)\quad(f\in L_\mathbb{R}^2(\mathfrak{N},dF))\label{eq3.22}
\end{equation}
where the kernel $K(x,t)$ corresponds to $\widetilde{\varphi^*\sigma\varphi}$ ($=\widetilde{K}$ \eqref{eq3.15}) in realization \eqref{eq3.21}. From $Y=-B-X$ one finds $\widetilde{Y}$,
\begin{equation}
(\widetilde{Y}f)(x)=a(x)f(x)+i\int\limits_\mathbb{R}\frac{K(x,t)}{t-x}dF(t)f(t)\quad(f\in L_\mathbb{R}^2(\mathfrak{N},dF))\label{eq3.23}
\end{equation}
where $a(x)=-b(x)-x$. Finally, using $A=YX$ one obtains
\begin{equation}
(\widetilde{A}f)(x)=xa(x)f(x)+i\int\limits_\mathbb{R}\frac{K(x,t)}{t-x}tdF(t)f(t)\quad(f\in L_\mathbb{R}^2(\mathfrak{N},dF)).\label{eq3.24}
\end{equation}
Equations \eqref{eq3.21} -- \eqref{eq3.24} are said to be the {\bf model representations} of the operators $X$, $B$, $Y$, $A$ written in terms of the spectral representation \eqref{eq3.21} of the operator $\widetilde{X}$.
\vspace{5mm}

{\bf 3.4.} Knowing self-adjoint operators $X$ and $Y$ from the factorization of $L(\lambda)$ \eqref{eq3.2}, one can give explicit formulas for solution to the Cauchy problem
\begin{equation}
\left\{
\begin{array}{lll}
\ddot{h}+B\dot{h}+Ah=f;\\
h(0)=h_0,\,\dot{h}(0)=h_1\,(0\leq t\leq T),
\end{array}\right.\label{eq3.25}
\end{equation}
here $h=h(.)$, $f=f(.)$ are vector functions from $H$ and $h_0$, $h_1$ are constant vectors from $H$. Let
\begin{equation}
X=\int\limits_\mathbb{R}\lambda dE_\lambda,\quad Y=\int\limits_\mathbb{R}wd\widetilde{E}_w\label{eq3.26}
\end{equation}
be spectral decompositions of self-adjoint operators $X$ and $Y$, and $\sigma(X)\cap\sigma(Y)=\emptyset$. The function
$$F_0(t)\stackrel{\rm def}{=}\int\limits_\mathbb{R}e^{\lambda t}dE_\lambda g+\int\limits_\mathbb{R}e^{wt}(wI-X)^{-1}d\widetilde{E}_w\widetilde{g},$$
where $g$, $\widetilde{g}\in H$, is correctly defined due to $\sigma(X)\cap\sigma(Y)=\emptyset$ and satisfies equation \eqref{eq3.25} since
$$\ddot{F}_0+BF_0+AF_0=\int\limits_\mathbb{R}L(\lambda)e^{\lambda t}dE_\lambda g+\int\limits_\mathbb{R}L(w)e^w(wI-X)^{-1}d\widetilde{E}_w\widetilde{g}$$
$$=\int\limits_\mathbb{R}(\lambda I-Y)(\lambda I-X)e^{\lambda t}dE_\lambda g+\int\limits_\mathbb{R}(wI-Y)e^{wt}d\widetilde{E}_w\widetilde{g}=0.$$
Initial data \eqref{eq3.25} give system of equations
$$\left\{
\begin{array}{lll}
g+K\widetilde{g}=h_0;\\
Xg+XK\widetilde{g}+\widetilde{g}=h_1
\end{array}\right.$$
where
\begin{equation}
K\stackrel{\rm def}{=}\int\limits_\mathbb{R}(wI-X)^{-1}d\widetilde{E}_w,\label{eq3.27}
\end{equation}
hence one finds that $\widetilde{g}=h_1-Xh_0$, $g=h_0-K(h_1-Xh_0)$. So, the function
\begin{equation}
F_0(t)=\int\limits_\mathbb{R}e^{\lambda t}dE_\lambda(h_0-K(h_1-Xh_0))+\int\limits_\mathbb{R}e^{wt}(wI-X)^{-1}d\widetilde{E}_w(h_1-Xh_0)\label{eq3.28}
\end{equation}
is the solution to Cauchy problem \eqref{eq3.25} when $f=0$.

Consider a function

$$F_1(t)=\int\limits_0^tds\int\limits_\mathbb{R}e^{\lambda(t-s)}dE_\lambda g(s)+\int\limits_0^tds\int\limits_\mathbb{R}e^{w(t-s)}(wI-X)d\widetilde{E}_w\widetilde{g}(s)$$
where $g(\cdot)$, $\widetilde{g}(\cdot)$ are vector functions from $H$. Obviously, $F_1(0)=0$.  Since
$$\dot{F}_1(t)=\int\limits_0^tds\int\limits_\mathbb{R}\lambda e^{\lambda-s}dE_\lambda g(s)+\int\limits_0^tds\int\limits_\mathbb{R}we^{w(t-s)}(wI-X)^{-1}d\widetilde{E}_w\widetilde{g}(s)$$
$$+g(t)+K\widetilde{g}(t),$$
then assuming that
\begin{equation}
g(t)+K\widetilde{g}(t)=0\label{eq3.29}
\end{equation}
one obtains that $\dot{F}_1(0)=0$ also. Using
$$\ddot{F}_1(t)=\int\limits_0^tds\int\limits_\mathbb{R}\lambda^2e^{\lambda(t-s)}dE_\lambda g(s)+\int\limits_0^tds\int\limits_\mathbb{R}w^2e^{w(t-s)}(wI-X)^{-1}d\widetilde{E}_w\widetilde{g}(s)$$
$$+X\dot{g}(t)+XK\widetilde{g}(t)+\widetilde{g}(t),$$
one obtains that the equation is equivalent to the equation
$$X(g+K\widetilde{g})+\widetilde{g}=f,$$
this, due to \eqref{eq3.29}, gives $\widetilde{g}=f$. So, the function
\begin{equation}
F_1(t)=-\int\limits_0^tds\int\limits_\mathbb{R}e^{\lambda(t-s)}dE_\lambda Kf(s)+\int\limits_0^t\int\limits_\mathbb{R}e^{w(t-s)}(wI-X)^{-1}d\widetilde{E}_wf(s)\label{eq3.30}
\end{equation}
is the solution to the non-homogenous ($f\not=0$) Cauchy problem \eqref{eq3.25} with zero initial data $h_0=h_1=0$.

\begin{theorem}\label{t3.2}
Let a self-adjoint operator $X$ be a spectral root of the pencil $L(\lambda)$ \eqref{eq3.2} and \eqref{eq3.26} be spectral decompositions of the operators $X$ and $Y$ from \eqref{eq3.2}, then solution to the Cauchy problem \eqref{eq3.25} is $h=F_0+F_1$ where $F_0$ and $F_1$ are given by \eqref{eq3.28} and \eqref{eq3.30} respectively, and the operator $K$ is \eqref{eq3.27}.
\end{theorem}
\vspace{5mm}

{\bf 3.5} Proceed to the second case (the {\bf elliptic condition}), when $C=0$, \eqref{eq3.8} takes place, and $X=iD$ ($D=D^*$).

First, consider a homogenous equation $\{D,B\}=0$.

\begin{remark}\label{r3.3}
If $\{D,B\}=0$, then
\begin{equation}
\begin{array}{lll}
({\rm i})\,D^2B=BD^2;\quad B^2D=DB^2;\\
({\rm ii})\,D\Ker B\subseteq\Ker B,\quad B\Ker D\subseteq \Ker D.
\end{array}\label{eq3.31}
\end{equation}
\end{remark}

Using spectral decomposition of $B$, one obtains
\begin{equation}
B=B_+-B_-\quad(B_\pm\stackrel{\rm def}{=}\pm\int\limits_{\mathbb{R}_\pm}\lambda dE_\lambda\geq0)\label{eq3.32}
\end{equation}
($E_\lambda$ is decomposition of the identity of the operator $B$) and
\begin{equation}
H=H_+\oplus H_-\oplus H_0\quad(H_\pm\stackrel{\rm def}{=}\overline{B_\pm H},\,H_0\stackrel{\rm def}{=}\Ker B).\label{eq3.33}
\end{equation}

\begin{lemma}\label{l3.3}
If $D$ is a bounded self-adjoint operator such that $\{D,B\}=0$, then $D_\pm\stackrel{\rm def}{=}P_\pm DP_\pm=0$, where $P_\pm$ are orthogonal projections onto $H_\pm$ of decomposition \eqref{eq3.33} corresponding to the operator $B$ \eqref{eq3.32}.
\end{lemma}

P r o o f. Equations $[P_+,B]=0$ and $[D,B]=0$ imply that $\{D_+,B_+\}=0$, therefore
$$0\leq\langle B_+D_th,D_+h\rangle=-\langle D_+^2B_+h,h\rangle\quad(\forall h\in H_+)$$
due to $B_+\geq0$. Since $B_+\geq0$ and $[D_+^2,B_+]=0$ \eqref{eq3.31}, then $D_+^2B_+=\sqrt{B_+}D_+^2\sqrt{B_+}\geq0$, and thus $\langle B_+D_+h,D_+h\rangle=0$, i. e., $B_+D_+h=0$, which gives $D_+h\in H_0$ and this is possible only when $D_+=0$. For $D_-$, the proof is analogous. $\blacksquare$

Taking into account (ii) \eqref{eq3.31} and Lemma \ref{l3.3}, write the operator $D$ block-by-block with regard to decomposition \eqref{eq3.33},
$$D=\left[
\begin{array}{ccc}
0&\Gamma&0\\
\Gamma^*&0&0\\
0&0&D_0
\end{array}\right]\quad(\Gamma=P_+DP_-,\,D_0=P_0DP_0),$$
then $\{D,B\}=0$ implies
\begin{equation}
B_+\Gamma-\Gamma B_-=0,\label{eq3.34}
\end{equation}
here $B_\pm$ is given by \eqref{eq3.32}. Since $\Gamma^*B_+-B_-\Gamma^*=0$, then $[B_+,\Gamma\Gamma^*]=0$ and $[B_-,\Gamma^*\Gamma]=0$, therefore \cite{26}
\begin{equation}
[B_-,|\Gamma|]=0,\quad[B_+,|\Gamma^*|]=0\label{eq3.35}
\end{equation}
where $|\Gamma|=\sqrt{\Gamma^*\Gamma}$, $|\Gamma^*|=\sqrt{\Gamma\Gamma^*}$. Polar representation of the operator $\Gamma$ is given by \cite{3, 26}
\begin{equation}
\Gamma=V|\Gamma|\quad(\Gamma^*=V^*|\Gamma^*|=|\Gamma|V^*)\label{eq3.36}
\end{equation}
where $V$ is a unitary operator from $\overline{|\Gamma|H}$ onto $\overline{|\Gamma^*|H}$. Equations \eqref{eq3.34} -- \eqref{eq3.36} imply that
$$0=B_+\Gamma-\Gamma B_-=(B_+V-VB_-)|\Gamma|,$$
therefore the restrictions $\left.B_-\right|_{\overline{|\Gamma|H}}$ and $\left.B_+\right|_{\overline{|\Gamma^*|H}}$ are unitarily equivalent. Instead of \eqref{eq3.33}, consider decomposition
\begin{equation}
H=G_+\oplus G_-\oplus G_0\quad(\dim G_+=\dim G_-)\label{eq3.37}
\end{equation}
where $G_+=\overline{|\Gamma^*|H_+}$, $G_-=\overline{|\Gamma|H_-}$, and
\begin{equation}
G_0=H_0\oplus\Ker|\Gamma|\oplus\Ker|\Gamma^*|=\{f+g:f\in\Ker B,g\in\Ker D\},\label{eq3.38}
\end{equation}
besides, $BG_0\subset G_0$, $DG_0\subset G_0$.

\begin{theorem}\label{t3.3}
Let $B$ and $D$ be a pair of anti-commuting, $\{D,B\}=0$, bounded self-adjoint operators in $H$. Then there is such decomposition \eqref{eq3.37}, \eqref{eq3.38} of a Hilbert space $H$ that

${\rm (i)}$ operators $B_-\geq0$, $|\Gamma|\geq0$ exist in $G_-$ ($\Ker B_-=\Ker|\Gamma|=0$) and $[B_-,|\Gamma|]=0$;

${\rm (ii)}$ there exists a unitary operator $V:$ $G_-\rightarrow G_+$ such that in terms of decomposition \eqref{eq3.37} operators $B$ and $D$ are given by
\begin{equation}
B=\left[
\begin{array}{ccc}
VB_-V^*&0&0\\
0&-B_-&0\\
0&0&B_0
\end{array}\right];\quad D=\left[
\begin{array}{ccc}
0&V|\Gamma|&0\\
|\Gamma|V^*&0&0\\
0&0&D_0
\end{array}\right],\label{eq3.39}
\end{equation}
besides, $\{B_0,D_0\}=0$.
\end{theorem}

\begin{remark}
Class of bounded self-adjoint anti-commuting operators $B$, $D$ is parameterized by the operator triple $\{B_-,|\Gamma|,V\}$ where $B_-\geq0$, $|\Gamma|\geq0$ ($[B_-,|\Gamma|]=0$) are operators in $G_-$, and $V:$ $G_-\rightarrow G_+$ is a unitary operator. The subspace $G_0$ \eqref{eq3.38} is characterized by $G_0=\Ker BD$ ($=\Ker DB$).
\end{remark}

\begin{corollary}\label{c3.1}
An arbitrary pair of anti-commuting self-adjoint bounded operators $B$, $D$ in terms of the decomposition \eqref{eq3.37} are unitarily equivalent to the operators
\begin{equation}
B=\left[
\begin{array}{ccc}
B_-&0&0\\
0&-B_-&0\\
0&0&B_0
\end{array}\right];\quad D=\left[
\begin{array}{ccc}
0&|\Gamma|&0\\
|\Gamma|&0&0\\
0&0&D_0
\end{array}\right]\label{eq3.40}
\end{equation}
where $G_0$ is given by \eqref{eq3.38} and $B_-$ and $|\Gamma|$ are non-negative commuting operators in $G_-$. Operators $B_0$ and $D_0$ are such that $B_0D_0=D_0B_0=0$.
\end{corollary}

Using Theorem \ref{t2.1}, realize the operator $B_-$ in $G_-$ by the operator of multiplication by independent variable ($G_-$ is isomorphic to $L_{\mathbb{R}_+}^2(\mathfrak{N},dF)$), besides, $|\Gamma|$ transforms (Theorem \ref{t2.2}) into the operator of multiplication \eqref{eq2.10} by the operator-valued function $n(x):$ $\mathfrak{N}\rightarrow{N}$, such that $n(x)\geq0$ ($\forall x$).

\begin{theorem}\label{t3.4}
Let $B$, $D$ be a pair of anti-commuting, $\{D,B\}=0$, bounded self-adjoint operators in $H$, besides, $\Ker DB=\{0\}$ ($G_0=\{0\}$ \eqref{eq3.38}). Then $H$ is unitarily equivalent to the space $L_{\mathbb{R}_+}^2(\mathfrak{N},dF)\otimes\mathbb{C}^2$, and $B$ and $D$, respectively, are unitarily equivalent to the operators
\begin{equation}
(\widetilde{B}f)(x)=\left[
\begin{array}{ccc}
I_\mathfrak{N}&0\\
0&-I_\mathfrak{N}
\end{array}\right]xf(x);\quad (\widetilde{D}f)(x)=\left[
\begin{array}{ccc}
0&n(x)\\
n(x)&0
\end{array}\right]f(x);\label{eq3.41}
\end{equation}
where $n(x)$ is an operator-valued function in $\mathfrak{N}$ such that $n(x)\geq0$ and $dFn=ndF$.
\end{theorem}

Thus the Pauli matrices
$$\sigma_1=\left[
\begin{array}{ccc}
0&1\\
1&0
\end{array}\right];\quad\sigma_3=\left[
\begin{array}{ccc}
1&0\\
0&-1
\end{array}\right]$$
are the primary (simplest) generators of anti-commuting self-adjoint operators. This result was obtained by Yu. S. Samoilenko \cite{30} when constructing muilti-parameter spectral decomposition of a pair of anti-commuting self-adjoint operators.

Functional realization of the block representation \eqref{eq3.41} is as follows. Let $F_s$ be the odd extension of $F$ onto $\mathbb{R}_-$, $F(-x)=-F(x)$ ($\forall x\in\mathbb{R}_+$), thus support of $dF_s$ is symmetric relative to zero. Matching the vector function $f(x)=\col[f_1(x),f_2(x)]\in L_{\mathbb{R}_+}^2(\mathfrak{N},dF)\otimes\mathbb{C}^2$ with the function $\widetilde{f}(x)=f_1(x)\chi_+(x)+f_2(-x)\chi_-(x)$ ($\chi_\pm(x)$ are characteristic functions of $\mathbb{R}_\pm$), one obtains the space $L_{\mathbb{R}_+}^2(\mathfrak{N},dF)\otimes\mathbb{C}^2$ isomorphic to the space $L_\mathbb{R}^2(\mathfrak{N},dF_s)$.

\begin{theorem}\label{t3.5}
If $B$, $D$ is a pair of anti-commuting bounded self-adjoint operators in $H$ and $\Ker DB=\{0\}$, then $H$ is isomorphic to the space $L_\mathbb{R}^2(\mathfrak{N},F_s)$, and $B$, $D$ are unitarily equivalent to the operators
\begin{equation}
(\widetilde{B}f)(x)=xf(x);\quad(\widetilde{D}f)(x)=n(x)f(-x)\label{eq3.42}
\end{equation}
where $F_s(x)$ is an odd function on $\mathbb{R}$, $n(x)$ is an operator-valued function in $\mathfrak{N}$, and $n(x)\geq0$, $n(-x)=n(x)$, $ndF_s=F_sn$.
\end{theorem}
\vspace{5mm}

{\bf 3.6.} Consider the non-homogenous equation $\{D,B\}=\varphi^*\sigma\varphi$ \eqref{eq3.8}. The subspace
\begin{equation}
H_1\stackrel{\rm def}{=}\span\{B^n\varphi^*\sigma\varphi h:h\in H,n\in\mathbb{Z}_+\}\label{eq3.43}
\end{equation}
is invariant relative to $B$, and $\varphi^*\sigma\varphi$ and let $H_2=H_1^\perp$, then in terms of decomposition $H=H_1\oplus H_2$ operators $B$ and $D$ are
$$B=\left[
\begin{array}{ccc}
B_1&0\\
0&B_2
\end{array}\right],\quad D=\left[
\begin{array}{ccc}
D_1&\Gamma\\
\Gamma^*&D_2
\end{array}\right]$$
where $B_k=P_kHP_k$, $D_k=P_kDP_k$ ($P_k$ are orthogonal projections onto $H_k$, $k=1$, 2), $\Gamma=P_1DP_2$. Equation $\{D,B\}=-\varphi^*\sigma\varphi$ yields
\begin{equation}
\left\{
\begin{array}{lll}
({\rm{i}})\,\{D_1,B_1\}=\varphi^*\sigma\varphi;\\
({\rm{ii}})\,B_1\Gamma+\Gamma B_2=0;\\
({\rm{iii}})\,\{D_2,B_2\}=0.
\end{array}\right.\label{eq3.44}
\end{equation}
Equation (iii) \eqref{eq3.44} is considered in Subsection 3.5. Turn to equation (iii) \eqref{eq3.44}.

\begin{remark}\label{r3.5}
Let $H_{1,0}=\Ker B_1$ and $H_{1,1}=H_1\ominus H_{1,0}$, then $H_1=H_{1,1}\oplus H_{1,0}$ and the operators $B_1$ and $D_1$ in this decomposition are
$$B_1=\left[
\begin{array}{ccc}
B_{1,1}&0\\
0&0
\end{array}\right],\quad D_1=\left[
\begin{array}{ccc}
D_{1,1}&\Gamma_1\\
\Gamma_1^*&D_{1,0}
\end{array}\right].$$
$\{D_1,B_1\}=-\varphi^*\sigma\varphi$ implies that $B_{1,1}\Gamma_1=-P_{1,1}\varphi^*\sigma\varphi P_{1,0}$. So, knowing $B_{1,1}$ as the solution to the equation $\{B_{1,1},D_{1,1}\}=-P_{1,1}\varphi^*\sigma\varphi P_{1,1}$ one can define the block $\Gamma_1$. Besides, $D_{1,0}$ is an arbitrary self-adjoint operator. Hereinafter, assume that $\Ker B_1=\{0\}.$
\end{remark}

\begin{theorem}\label{t3.6}
Let $H=H_1$ \eqref{eq3.43} be separable, then for all the bounded self-adjoint operators $B$ and $\varphi^*\sigma\varphi$ such that $B\geq0$, $\Ker B=\{0\}$  and $\varphi^*\sigma\varphi$ is $B$-bounded, there exists the unique self-adjoint bounded solution $D$ to the equation $\{D,B\}=-\varphi^*\sigma\varphi$.

In the spectral representation (Theorem \ref{t2.1}) of the operator $\widetilde{B}$ ($=\widetilde{T}$ \eqref{eq3.44}) in $L_{\Omega}^2(\mathfrak{N},dF)$ ($\Omega=\supp dF\subset\mathbb{R}_+$ is bounded), the only self-adjoint bounded solution $\widetilde{D}$ to the equation $\{\widetilde{D},\widetilde{B}\}=-\widetilde{\varphi^*\sigma\varphi}$ ($\widetilde{\varphi^*\sigma\varphi}=\widetilde{K}$ \eqref{eq3.15}) is given by
\begin{equation}
(\widetilde{D}f)(x)=-\int\limits_{\mathbb{R}_+}\frac{K(x,t)}{t+x}dF(t)f(t)\label{eq3.45}
\end{equation}
where the operator-valued function $K(x,t)$ in $\mathfrak{N}$ has the same properties as in Lemma \eqref{l3.2}.
\end{theorem}

P r o o f. Boundedness of $\widetilde{D}$ \eqref{eq3.45} follows from $B$-boundedness of the operator $\varphi^*\sigma\varphi$, moreover, $\widetilde{D}$ is self-adjoint, due to $K^*(x,t)=K(t,x)$, and satisfies the relation $\{\widetilde{D},\widetilde{B}\}=-\widetilde{\varphi^*\sigma\varphi}$. Uniqueness of the solution $\widetilde{D}$ follows from Lemma \ref{l3.3}. Really, if $\widetilde{D}_1$ is another solution to the equation $\{\widetilde{D},\widetilde{B}\}=-\widetilde{\varphi^*\sigma\varphi}$, then $\{(D-\widetilde{D}),B\}=0$, which in view of $\widetilde{B}\geq0$ and of $\Ker\widetilde{B}=0$ gives $\widetilde{D}=\widetilde{D}_1$. $\blacksquare$

\begin{remark}\label{r3.6}
In proper formulation, Theorem \ref{t3.6} takes place for $B\leq0$, $\Ker B=\{0\}$ also. Operator $\widetilde{D}$ \eqref{eq3.45} is an analogue of the well-known \cite{21} {\bf Stieltjes transform}.
\end{remark}

Consider general case, without supposing that $B\geq0$ ($\leq0$), and let $H=H_1$ \eqref{eq3.43} and $\Ker B=\{0\}$, and $\widetilde{B}$ ($=\widetilde{T}$ \eqref{eq3.14}) be spectral realization of the operator $B$ in the space $L_\Omega^2(\mathfrak{N},dF)$ where $\Omega=\supp dF\subset\mathbb{R}$ is bounded, besides, $\widetilde{\varphi^*\sigma\varphi}$ ($=\widetilde{K}$ \eqref{eq3.15}) is the corresponding representation in this realization. General solution $\widetilde{D}_1$ to the equation $\{\widetilde{D}_1,\widetilde{B}\}=-\widetilde{\varphi^*\sigma\varphi}$ equals $\widetilde{D}_1=\widetilde{D}_1+N$ (see proof of Theorem \ref{t3.6}) where $\widetilde{D}$ is given by \eqref{eq3.45} and $N$ is a self-adjoint operator, anti-commuting with $\widetilde{B}$, $\{N,\widetilde{B}\}=0$.

Existence of an anti-commuting with $\widetilde{B}$ operator $N$ means (Theorem \ref{t3.5}) that support $\Omega$ of measure $dF$ has symmetric part $\Omega_s\subseteq\Omega$,
\begin{equation}
\Omega_s=\Omega_-\cup\Omega_+\quad(\Omega_\pm\subset(\Omega\cap\mathbb{R}_\pm),\,x\in\Omega_+\Leftrightarrow-x\in\Omega_-).\label{eq3.46}
\end{equation}
Block structure of $\widetilde{B}$, $N$, in view of \eqref{eq3.39} ($\Ker\widetilde{B}=\{0\}$), in the space $L_{\Omega_s}^2(\mathfrak{N},dF)=L_{\Omega_-}^2(\mathfrak{N},dF)\oplus L_{\Omega_+}^2(\mathfrak{N},dF)$ is
$$\widetilde{B}_s=\left[
\begin{array}{ccc}
V\widetilde{B}_-V^*&0\\
0&-\widetilde{B}_-
\end{array}\right];\quad N=\left[
\begin{array}{ccc}
0&V|\Gamma|\\
|\Gamma|V^*&0
\end{array}\right]$$
where commuting operators $\widetilde{B}_-$ and $|\Gamma|$ in $L_{\Omega_-}^2(\mathfrak{N},dF)$ are
$$(\widetilde{B}_-f_-)(x)=-xf_-(x);\quad(|\Gamma|f_-)(x)=n(x)f_-(x)$$
($\forall f_-\in L_{\Omega_-}^2(\mathfrak{N},dF)$), besides, $n(x)$ is an operator-valued function in $\mathfrak{N}$ such that $n(x)\geq0$, $n(x)dF(x)=dF(x)n(x)$ ($\forall x\in\Omega_-$). Unitary operator $V:$ $L_{\Omega_-}^2(\mathfrak{N},dF)\rightarrow L_{\Omega_+}^2(\mathfrak{N},dF)$ is given by
$$(Vf_-)(x)=v(x)f_-(-x)\quad(x\in\Omega_+,\,f_-\in L_{\Omega_-}^2(\mathfrak{N},dF)$$
where $v(x)$ is an invertible ($\forall x\in\Omega_+$) operator-valued function in $\mathfrak{N}$. This representation for $V$ follows from the permutability of $VS$ and $B$ where $Sf_-(x)=f_-(-x)$ ($L_{\Omega_-}^2(\mathfrak{N},dF)\rightarrow L_{\Omega_+}^2(\mathfrak{N},dF)$) due to $VB_-V^*=B_-$. Unitarity of $V$ means that
\begin{equation}
v^*(x)dF(x)v(x)=dF(-x)\quad(\forall x\in\Omega_+),\label{eq3.47}
\end{equation}
besides, $V^*$: $L_{\Omega_+}^2(\mathfrak{N},dF)\rightarrow L_{\Omega_-}^2(\mathfrak{N},dF)$ is
$$(V^*f_+)(x)=v^{-1}(-x)f_+(-x)\quad(x\in\Omega_-,f_+\in L_{\Omega_+}^2(\mathfrak{N},dF)).$$
If $\Phi=\col[f_=,f_-]$ ($f_\pm\in L_{\Omega_\pm}^2(\mathfrak{N},dF)$), then
$$(\widetilde{B}_s\Phi)(x)=x\Phi(x),\quad(N\Phi)(x)=\col[v(x)n(-x)f_-(x),n(x)v^{-1}(-x)f_+(-x)].$$

\begin{theorem}\label{t3.7}
If $H$ is separable and $H=H_1$ \eqref{eq3.43}, then for all the bounded self-adjoint operators $B$ and $\varphi^*\sigma\varphi$ such that $\Ker B=\{0\}$ and $\varphi^*\sigma\varphi$ is $B$-bounded there exists the bounded self-adjoint solution $D$ to the equation $\{D,B\}=-\varphi^*\sigma\varphi$. In the spectral representation (Theorem \ref{t2.1}) of the operator $\widetilde{B}$ ($=\widetilde{T}$ \eqref{eq3.14}) in $L_\Omega^2(\mathfrak{N},dF)$ ($\Omega=\supp dF$) the solution $\widetilde{D}$ to the equation $\{\widetilde{D},\widetilde{B}\}=-\widetilde{\varphi^*\sigma\varphi}=\widetilde{K}$ \eqref{eq3.15}) is
\begin{equation}
(\widetilde{D}f)(x)=m(x)f_s(-x)-\int\limits_\Omega\frac{K(x,t)}{t+x}dF(t)f(t)\label{eq3.48}
\end{equation}
where $f_s=f\chi_{\Omega_s}$ ($\chi_{\Omega_s}$ is characteristic function of the set $\Omega_s$ \eqref{eq3.46}, $\Omega_s\subseteq\Omega$), $m(x)$ is an operator-valued function in $\mathfrak{N}$ given on $\Omega_s$,
\begin{equation}
m(x)\stackrel{\rm def}{=}\left\{
\begin{array}{ccc}
n(x)v^{-1}(-x)&(x\in\Omega_-);\\
v(x)n(-x)&(x\in\Omega_+),
\end{array}\right.\label{eq3.49}
\end{equation}
besides, $n(x)\geq0$ is an operator-valued function in $\mathfrak{N}$ on $\Omega_-$ and $dFn=ndF$; $v(x)$ is an invertible operator-valued function in $\mathfrak{N}$ given on $\Omega_+$ for which \eqref{eq3.47} takes place; $K(x,t)$ is an operator-valued function in $\mathfrak{N}$ with the same property as in Lemma \ref{l3.2}.
\end{theorem}

\begin{remark}\label{r3.7}
If $v(x)=I_\mathfrak{N}$ ($x\in\Omega_+$), then measure $dF$ is even on $\Omega_s$ due to \eqref{eq3.47} and $m(x)$ coincides with even extension of $n(x)$ onto $\Omega_+$ due to \eqref{eq3.49}.

In the case of $n(x)=I_\mathfrak{N}$, the function $m(x)$ is defined by the function $v(x)$ \eqref{eq3.49} and \eqref{eq3.47} takes place.

Non-integral summand in \eqref{eq3.48} vanishes if symmetric subspace $\Omega_s$ is absent in $\Omega$ or when there is no invertible operator-valued function $v(x)$ for which \eqref{eq3.47} takes place.
\end{remark}
\vspace{5mm}

{\bf 3.7} Proceed to the relation (ii) \eqref{eq3.44}, $B_1\Gamma+\Gamma B_2=0$. Using polar decomposition \eqref{eq3.36} for $\Gamma$ and self-adjointness of $B_1$, $B_2$, obtain
\begin{equation}
[B_1,|\Gamma^*|]=0,\quad[B_2,|\Gamma|]=0,\label{eq3.50}
\end{equation}
therefore
$$\{B_1V+VB_2\}|\Gamma|=0$$
where $V$ is a unitary operator from $\overline{|\Gamma|H_2}$ onto $\overline{|\Gamma^*|H_1}$. Thus, $B_2f=-V^*B_1Vf$ ($\forall f\in\overline{|\Gamma|H_2}$). Equation \eqref{eq3.50} implies $B_1\Ker|\Gamma^*|\subset\Ker|\Gamma^*|$ and $B_2\Ker|\Gamma|\subset\Ker|\Gamma|$.

Summarizing results of Subsections 3.5, 3.6, gives the statement.

\begin{theorem}\label{t3.8}
Let $H$ be separable and $B$ and $\varphi^*\sigma\varphi$ be bounded self-adjoint operators ($\varphi^*\sigma\varphi$ be $B$-bounded). Then there exists the decomposition
\begin{equation}
H=H_1\oplus H_2\oplus H_0\label{eq3.51}
\end{equation}
($H_1$ is given by \eqref{eq3.43}) such that self-adjoint bounded operators satisfying the relation $\{D,B\}=-\varphi^*\sigma\varphi$ have the block structure
\begin{equation}
B=\left[
\begin{array}{ccc}
B_1&0&0\\
0&-V^*B_1V&0\\
0&0&B_0
\end{array}\right];\quad D=\left[
\begin{array}{ccc}
D_1&|\Gamma^*|V&0\\
V^*|\Gamma^*|&D_2&0\\
0&0&D_0
\end{array}\right].\label{eq3.52}
\end{equation}
Besides,

${\rm(i)}$ self-adjoint operators $B_1$, $D_1$, $|\Gamma^*|$ ($\geq0$) in $H_1$ have the properties $\{D_1,B_1\}=-\varphi^*\sigma_1\varphi$ (description of such pairs is given in Theorem \ref{t3.7}) and $[B_1,|\Gamma^*|]=0$ (model representation of such operators is obtained in Theorem \ref{t2.2});

${\rm(ii)}$ there are a unitary operator $V$ from $H_2$ onto $\overline{|\Gamma^*|H_1}$ and a self-adjoint bounded operator $D_2$ in $H_2$ such that $\{D_2,V^*B_1V\}=0$ (model realization of these operators is obtained in Theorem \ref{t3.4});

${\rm(iii)}$ $B_0$ and $D_0$ are anti-commuting, $\{D_0,B_0\}=0$, self-adjoint operators in $H_0$, structure of which is given in Theorem \ref{t3.3}.
\end{theorem}

\begin{remark}\label{r3.8}
Under the suppositions of Theorem \ref{t3.8}, the roots of the quadratic pencil are
\begin{equation}
X=iD,\quad Y=-B-iD\label{eq3.53}
\end{equation}
where $B$ and $D$ are bounded self-adjoint operators. If $B\geq\delta I$ ($\delta>0$), then $\sigma(Y)$ lies in the half-plane $\Re z\leq\-\delta$ and thus, in this case, the root $X$ is a spectral one. Spectrality of the root $X$ also takes place if $B\leq-\delta I$ ($\delta>0$).
\end{remark}

\section{Characteristic function}

{\bf 4.1} Consider model representations \eqref{eq3.21} -- \eqref{eq3.24} of the operators $X$, $B$, $Y$, $A$ and suppose that
\begin{equation}
\dim\mathfrak{N}=1;\quad\supp dF=\Omega=[a,b]\,(-\infty<a<b<\infty),\quad dF(x)=dx,\label{eq4.1}
\end{equation}
then $L_\Omega^2(\mathfrak{N},dF)=L^2(a,b)$. Hence it follows that $\widetilde{X}$ \eqref{eq3.21} is a bounded self-adjoint operator with simple absolutely continuous spectrum and $\sigma(\widetilde{X})=[a,b]$. Let $\rank(\varphi^*\sigma\varphi)=1$, then $K(x,t)=Jv(x)\overline{v}(t)$ ($\sigma=J=\pm1$, $v\in L^2(a,b)$), besides,
\begin{equation}
(\widetilde{X}f)(x)=xf(x);\quad(\widetilde{Y}f)(x)=-a(x)f(x)+iJv(x)/\hspace{-4.3mm}\int\limits_a^b\frac{f(t)\overline{v}(t)}{t-x}dt\label{eq4.2}
\end{equation}
where $f\in L^2(a,b)$; $a(x)=x+b(x)$; $b(x)$ is a real function such that $\sup\limits_x|b(x)|<\infty$. The operator $\widetilde{\varphi}:$ $L^2(a,b)\rightarrow\mathbb{C}$ is
\begin{equation}
\widetilde{\varphi}f=\int\limits_a^bf(x)dx;\quad(\varphi^*\xi=\xi v(x),\,\forall\xi\in\mathbb{C}).\label{eq4.3}
\end{equation}
For the model \eqref{eq4.1} -- \eqref{eq4.3}, calculate the characteristic function $S_\Delta(\lambda,B)$ \eqref{eq1.26},
\begin{equation}
S_\Delta(\lambda,B)=I-i\widetilde{\varphi}\widetilde{L}^{-1}(\lambda)\widetilde{\varphi}^*\sigma=I-i\widetilde{\varphi}(\lambda I-\widetilde{X})^{-1}(\lambda I-\widetilde{Y})^{-1}\widetilde{\varphi}^*J\label{eq4.4}
\end{equation}
where $\widetilde{X}$, $\widetilde{Y}$ are roots \eqref{eq4.2} of the pencil $\widetilde{L}(\lambda)=\lambda^2I+\lambda\widetilde{B}+\widetilde{A}$ \eqref{eq3.2}. Let
\begin{equation}
f=(\lambda I-\widetilde{Y})^{-1}\widetilde{\varphi}J\xi,\label{eq4.5}
\end{equation}
then \eqref{eq4.2}, \eqref{eq4.3} imply
\begin{equation}
\zeta(\lambda,x)f(x)+Jv(x)\frac1i/\hspace{-4.3mm}\int\limits_a^b\frac{f(t)\overline{v}(t)}{t-x}dt=\zeta Jv(x)\label{eq4.6}
\end{equation}
where
\begin{equation}
\zeta(\lambda,x)\stackrel{\rm def}{=}\lambda+a(x)=\lambda+x+b(x).\label{eq4.7}
\end{equation}
Multiplying \eqref{eq4.6} by $\overline{v(x)}$ and denoting
\begin{equation}
F(x)\stackrel{\rm def}{=}f(x)v(x);\quad\omega(x)\stackrel{\rm def}{=}J|v(x)|^2,\label{eq4.8}
\end{equation}
one obtains
\begin{equation}
\zeta(\lambda,x)F(x)+\frac{\omega(x)}i/\hspace{-4.3mm}\int\limits_a^b\frac{F(t)dt}{t-x}=\zeta\omega(x).\label{eq4.9}
\end{equation}
Continue the functions $f$, $v$, $b$ (and thus $F$ and $\omega$ also) by zero outside the interval $[a,b]$. Define the function
\begin{equation}
\Phi(z)\stackrel{\rm def}{=}\frac1{2\pi i}\int\limits_\mathbb{R}\frac{F(t)}{t-z}dt\quad(z\in\mathbb{C})\label{eq4.10}
\end{equation}
where $F$ is given by \eqref{eq4.8}. Sokhotski formula \cite{20,21,22} for the Cauchy type integral \eqref{eq4.10} are true if $F\in L^p(\mathbb{R})$ ($p>1$). If
\begin{equation}
\sup\limits_x|v(x)|<\infty,\label{eq4.11}
\end{equation}
then $F\in L^2(\mathbb{R})$ (or $f\in L^2(\mathbb{R})$), and thus the Sokhotski formula \cite{22} are true:
\begin{equation}
\Phi_+(x)-\Phi_-(x)=F(x);\quad \Phi_+(x)+\Phi_-(x)=\frac1{\pi i}/\hspace{-4.3mm}\int\limits_\mathbb{R}\frac{F(t)dt}{t-x}\label{eq4.12}
\end{equation}
where $\Phi_\pm(x)=\Phi(x\pm i0)$ are boundary values on $\mathbb{R}$ from $\mathbb{C}_\pm$ of the function $\Phi(z)$ \eqref{eq4.10}. Substituting \eqref{eq4.12} into \eqref{eq4.9} gives
$$(\zeta(\lambda,x)+\pi\omega(x))\Phi_+(x)=(\zeta(\lambda,x)-\pi\omega(x))\Phi_-(x)+\xi\omega(x)$$
and one arrives at the Riemann boundary problem \cite{20,21}
\begin{equation}
\Phi_+(x)=d(x,\lambda)\Phi_-(x)+\frac{\xi\omega(x)}{\zeta(x,\lambda)+\pi\omega(x)},\label{eq4.13}
\end{equation}
here
\begin{equation}
d(x,\lambda)\stackrel{\rm def}{=}\frac{\zeta(\lambda,x)-\pi\omega(x)}{\zeta(\lambda,x)+\pi\omega(x)},\label{eq4.14}
\end{equation}
besides, $d(x,\lambda)=1$ $\forall x\not\in[a,b]$ and $\forall\lambda\not=0$ from $\mathbb{C}$. For each $\lambda$, such that $\zeta(\lambda,x)\pm\pi\omega(x)\not=0$ ($\forall x\in\mathbb{R}$) determine the canonic function $X(z,\lambda)$ \cite{20,21} as the solution of the homogenous boundary problem
\begin{equation}
X_+(x,\lambda)=d(x,\lambda)X_-(x,\lambda)\label{eq4.15}
\end{equation}
normalized by the condition $\chi(\infty,\lambda)=1$. Taking the logarithm of \eqref{eq4.15} and using the Sokhotski formula gives \cite{20,21}
\begin{equation}
X(z,\lambda)=\exp\{\Gamma(z,\lambda)\};\quad\Gamma(z,\lambda)\stackrel{\rm def}{=}\frac1{2\pi i}\int\limits_\mathbb{R}\frac{\ln d(x,\lambda)}{x-z}dx.\label{eq4.16}
\end{equation}
Equations \eqref{eq4.13}, \eqref{eq4.15} imply
\begin{equation}
\frac{\Phi_+(x)}{X_+(x,\lambda)}=\frac{\Phi_-(x)}{X_-(x,\lambda)}+\frac{\zeta\omega(x)}{X_+(x,\lambda)(\zeta(x,\lambda)+\pi\omega(x))}\label{eq4.17}
\end{equation}
Set the function
\begin{equation}
\psi(z,\lambda)\stackrel{\rm def}{=}\frac1{2\pi i}\int\limits_\mathbb{R}\frac{\zeta\omega(x)}{\chi_+(x,\lambda)(\zeta(x,\lambda)+\pi\omega(x))}\frac{dx}{x-z}\quad(z\in\mathbb{C})\label{eq4.18}
\end{equation}
where $\omega\in L^2(\mathbb{R})$, due to \eqref{eq4.11}. Using Sokhotski formula for $\psi(x,\lambda)$, rewrite relation \eqref{eq4.17} as
$$\frac{\Phi_+(x)}{\chi_+(x,\lambda)}-\psi_+(x,\lambda)=\frac{\Phi_-(x)}{\chi_-(x,\lambda)}-\psi_-(x,\lambda).$$
Hence it follows that for every fixed $\lambda$ function $\Phi(z)/\chi(z,\lambda)-\psi(z,\lambda)$ is analytic for all $z\in\mathbb{C}$ and coincides with the polynomial $P_{\varkappa}(z,\lambda)$ (in $x$) \cite{20,21} of degree ${\varkappa}\geq0$, where ${\varkappa}$ is {\bf index of the Riemann problem} \eqref{eq4.13}. Therefore
\begin{equation}
P(z)=\chi(z,\lambda)\{\psi(z,\lambda)+P_{\varkappa}(z,\lambda)\}.\label{eq4.19}
\end{equation}
Expression \eqref{eq4.4} for $S_\Delta(\lambda,B)$, due to \eqref{eq4.5}, \eqref{eq4.3}, \eqref{eq4.10}, is
$$S_\Delta(\lambda,B)\zeta=\xi-i\widetilde{\varphi}(\lambda I-\widetilde{X})^{-1}f=\xi+i\int\limits_a^b\frac{f(t)\overline{v(t)}}{t-\lambda}dt=\xi-2\pi\Phi(\lambda)$$
and taking into account \eqref{eq4.16}, \eqref{eq4.18}, \eqref{eq4.19} one has
$$S_\Delta(\lambda,B)\xi=\xi-2\pi\exp\{\Gamma(\lambda,\lambda)\}\cdot\left\{\frac1{2\pi i}\int\limits_\mathbb{R}\frac{\xi\omega(x)}{X_+(x,\lambda)(\zeta(x,\lambda)+\pi\omega(x))}\right.$$
$$\left.\times\frac{dx}{x-\lambda}+P_\varkappa(\lambda,\lambda)\right\}.$$
Since
$$X_+^{-1}(x,\lambda)=\exp\{-\Gamma(x,\lambda)\}=[d(x,\lambda)]^{-1/2}\exp\left\{-/\hspace{-4.3mm}\int\limits_\mathbb{R}\frac{\ln d(t,\lambda)}{t-x}dt\right\},$$
then
$$S_\Delta(\lambda,B)\xi=\xi\left\{1+i\int\limits_\mathbb{R}\frac{\omega(x)}{[\zeta^2(x,\lambda)-\pi^2\omega^2(x)]^{1/2}}\cdot\exp\left\{\frac1{2\pi i}\int\limits_\mathbb{R}\frac{\ln d(t,\lambda)dt}{t-\lambda}\right.\right.$$
$$\left.\left.-\frac1{2\pi i}/\hspace{-4.3mm}\int\limits_\mathbb{R}\frac{\ln(t,\lambda)}{t-x}dt\right\}\cdot\frac{dx}{x-\lambda}-2\pi P_\varkappa(\lambda,\lambda)\exp\{\Gamma(\lambda,\lambda)\}\right\}.$$
Taking into account that
$$S_\Delta(\lambda,B)=I-i\frac1{\lambda^2}\varphi^*\sigma\varphi+o\left(\frac1{\lambda^2}\right)\quad(|\lambda|\gg1),$$
one obtains that $P_2(\lambda,\lambda)=0$.

\begin{theorem}\label{t4.1}
Suppose that the spectral root $X$ of a quadratic pencil $L(\lambda)$ \eqref{eq3.2} is a bounded self-adjoint operator with simple absolutely continuous spectrum and $\sigma(X)=[a,b]$ ($-\infty<a<b<\infty$). And let $\rank(\varphi^*\sigma\varphi)=1$ and thus kernel of the operator $\widetilde{\varphi^*\sigma\varpi}$ \eqref{eq3.15} in spectral decomposition of $\widetilde{X}$ \eqref{eq4.2} be ($J=\pm1$), and \eqref{eq4.11} take place for $v$.

Then characteristic function $S_\Delta(\lambda,B)$ \eqref{eq4.4} of the model pair $\widetilde{X}$, $\widetilde{Y}$ \eqref{eq4.2} is
\begin{equation}
S_\Delta(\lambda,B)=1+i\int\limits_a^b\frac{\omega(x)}{[\xi^2(x,\lambda)-\pi^2\omega^2(x)]^{1/2}}\exp\left\{\frac{\lambda-x}{2\pi i}/\hspace{-4.3mm}\int\limits_\mathbb{R}\frac{\ln d(t,\lambda)}{(t-\lambda)(t-x)}\right\}\frac{dx}{x-\lambda}\label{eq4.20}
\end{equation}
where $\zeta(x,\lambda)$, $\omega(x)$, and $d(x,\lambda)$ are given by \eqref{eq4.7}, \eqref{eq4.8} and \eqref{eq4.14}, respectively.
\end{theorem}

As is known \cite{25}, a self-adjoint operator with simple absolutely continuous spectrum is realized by the operator of multiplication by independent variable in the space $L^2(\mathbb{R},\sigma'(x)dx)$. Condition $\sigma(X)=[a,b]$ implies that $\sigma'(x)>0$ ($\forall x\in[a,b]$) and $\sigma'(x)=0$ for all $x\not\in[a,b]$. Mapping $f\rightarrow f\sqrt{\sigma'}$ sets unitary isomorphism between the spaces $L^2(\mathbb{R},\sigma'(x)dx)$ and $L^2(a,b)$. Therefore, operator of multiplication by independent variable in $L^2[a,b)$ is a model of a self-adjoint operator with simple spectrum such that $\sigma(x)=[a,b]$.

\begin{remark}\label{r4.1}
Singularities of characteristic function $S_\Delta(\lambda,B)$ \eqref{eq4.20} (spectrum of the pencil $L(\lambda)$ \eqref{eq3.2}) is union $\sigma(X)\cup\sigma(Y)$, where $\sigma(X)=[a,b]$ and $\sigma(Y)$ coincides with the set of values of the function $\pm\pi(J|v(x)|^2)^{1/2}-x=b(x)$ when $x\in[a,b]$ which for $J=1$ are real and for $J=-1$ are complex. Functions $b$ and $v$ are such that $\sigma(X)\cap\sigma(Y)=\emptyset$, due to spectrality of the root $\widetilde{X}$.
\end{remark}
\vspace{5mm}

{\bf 4.2} Let conditions \eqref{eq4.1} hold and $\rank(\varphi^*\sigma\varphi)=r<\varphi$ where $\sigma=J$ ($=J^*=J^{-1}$) is an involution. Then the operators $\widetilde{X}$ \eqref{eq3.20} and $\widetilde{Y}$ \eqref{eq3.23} are
\begin{equation}
(\widetilde{X}f)(x)=xf(x);\quad(\widetilde{Y}f)(x)=-a(x)f(x)+i/\hspace{-4.3mm}\int\limits_a^b\frac{K(x,t)}{t-x}f(t)dt\label{eq4.21}
\end{equation}
where $f\in L^2(a,b)$; $a(x)=x+b(x)$ and $b(x)$ is a real function such that $\sup\limits_x|b(x)|<\infty$; the kernel $K(x,t)$ is
\begin{equation}
K(x,t)=\sum\limits_{\alpha,\beta=1}^rv_\alpha(x)j_{\alpha\beta}\overline{v_\beta(t)};\label{eq4.22}
\end{equation}
here $\{v_\alpha\}_1^r$ is a set of linearly independent functions in $L^2(a,b)$; $j_{\alpha\beta}$ are matrix elements, $j_{\alpha\beta}=\langle Je_\beta,e_\alpha\rangle$ of the involution $J$ in the orthonormal basis $\{e_\alpha\}_1^r$ of space $E$. Operator $\widetilde{\varphi}:$ $L^2(a,b)\rightarrow E$ is given by
\begin{equation}
\widetilde{\varphi}f=\int\limits_a^bf(t)\sum\limits_\beta\overline{v_\beta(t)}dte_\beta\quad(\varphi^*e_\beta=e_\beta,\,1\leq\beta\leq r)\label{eq4.23}
\end{equation}
and
$$\widetilde{\varphi}^*J\widetilde{\varphi}f=\int\limits_a^bK(x,t)f(t)dt$$
where the kernel $K(x,t)$ coincides with \eqref{eq4.22}.

Calculate the characteristic matrix function $S_\Delta(\lambda,B)=[\langle S_\Delta(\lambda,B)e_\alpha,e_\beta\rangle]$ ($S_\Delta(\lambda,B)$ is given by formula \eqref{eq4.4}). Let
\begin{equation}
f_\alpha\stackrel{\rm def}{=}(\lambda I-\widetilde{Y})\widetilde{\varphi}^*Je_\alpha\quad(1\leq\alpha\leq r),\label{eq4.24}
\end{equation}
then
$$\zeta(x,\lambda)f_\alpha(x)+\sum\limits_{\gamma,\beta}v_\gamma(x)J_{\gamma,\beta}\frac1i/\hspace{-4.3mm}\int\limits_a^b\frac{f_\alpha(t)\overline{v}_\beta(t)}
{t-x}=\sum\limits_{\gamma=1}^rJ_{\gamma,\beta}v_\gamma(x)$$
where $\zeta(x,\lambda)=\lambda+a(x)$. Multiply this equation by $\overline{v_\gamma(x)}$; then in terms of the functions
\begin{equation}
F_{\delta,\alpha}(x)=f_\alpha(x)\overline{v_\delta}(x);\quad\omega_{\delta,\gamma}(x)=v_\gamma(x)\overline{v_\delta(x)}\label{eq4.25}
\end{equation}
one obtains the system of equations
\begin{equation}
\zeta(x,\lambda)F_{\delta,\alpha}(x)+\sum\limits_{\gamma,\beta=1}^r\omega_{\delta,\gamma}(x)j_{\gamma,\beta}\frac1i/\hspace{-4.3mm}\int\limits_a^b\frac{F_{\beta,
\alpha}(t)}{t-x}dt=\sum\limits_\gamma\omega_{\delta,\gamma}(x)j_{\gamma,\alpha}.\label{eq4.26}
\end{equation}
Continue the functions $f_\alpha(\cdot)$, $v_\alpha(\cdot)$, $b(\cdot)$ (thus, $F_{t,\alpha}(\cdot)$, $\omega_{\delta,\gamma}(\cdot)$ also) with zero outside $[a,b]$ and define the matrices
\begin{equation}
F(x)\stackrel{\rm def}{=}\left[
\begin{array}{ccc}
F_{1,1}(x)&...&F_{1,r}(x)\\
...&...&...\\
F_{r,1}(x)&...&F_{r,r}(x)
\end{array}\right];\quad\Omega(x)\stackrel{\rm def}{=}\left[
\begin{array}{ccc}
\omega_{1,1}(x)&...&\omega_{1,r}(x)\\
...&...&...\\
\omega_{r,1}(x)&...&\omega_{r,r}(x)
\end{array}\right];\quad J=[j_{\beta,\alpha}],\label{eq4.27}
\end{equation}
then equations \eqref{eq4.26} in the matrix form become
\begin{equation}
\zeta(x,\lambda)F(x)+\Omega(x)J\frac1i/\hspace{-4.3mm}\int\limits_{\mathbb{R}}\frac{F(t)}{t-x}dt=\Omega(x)J.\label{eq4.28}
\end{equation}
Consider the matrix function
\begin{equation}
\Phi(z)=\frac1{2\pi i}\int\limits_{\mathbb{R}}\frac{F(t)dt}{t-z}\quad(z\in\mathbb{C}).\label{eq4.29}
\end{equation}
If
\begin{equation}
\sup\limits_{x}|v_\alpha(x)|<\infty\quad(1\leq\alpha\leq r),\label{eq4.30}
\end{equation}
then $F_{\delta,\alpha}L^2(\mathbb{R})$, and one can use Sokhotski formulas \cite{20,21} which are easily transferred to the matrix case,
\begin{equation}
\Phi_+(x)-\Phi_-(x)=F(x);\quad\Phi_+(x)-\Phi_-(x)=\frac{\pi i}/\hspace{-4.3mm}\int\limits_{\mathbb{R}}\frac{F(t)dt}{t-x}\label{eq4.31}
\end{equation}
where $\Phi_\pm(x)=\Phi(x\pm i0)$ are boundary values on $\mathbb{R}$ from $\mathbb{C}_\pm$ of the function $\Phi(z)$ \eqref{eq4.29}. Substituting \eqref{eq4.31} into \eqref{eq4.28} gives
\begin{equation}
\{\zeta(x,\lambda)I+\delta\Omega(x)J\}\Phi_+(x)=\{\zeta(x,\lambda)I-\pi\Omega(x)J\}\Phi_-(x)+\Omega(x)J.\label{eq4.32}
\end{equation}
Definition of $\Omega(x)$ \eqref{eq4.27} implies that
\begin{equation}
(\Omega(x)J)^2=p(x)\Omega(x)J\quad(p(x)=\sum\limits_{\alpha,\beta}v_\alpha(x)j_{\alpha,\beta}\overline{v_\beta(x)}),\label{eq4.33}
\end{equation}
therefore, for $\lambda\in\mathbb{C}\backslash\mathbb{R}$ and $|\lambda|\gg1$, one obtains
$$(\zeta I+\pi\Omega J)^{-1}=\frac1\zeta\left(I+\frac{\pi\Omega J}\zeta\right)^{-1}=\frac1\zeta\left(I-\frac{\pi\Omega\zeta}\zeta+\left(\frac{\pi\Omega J}\zeta\right)^2-...\right)$$
$$=\frac1\zeta I-\frac{\pi\Omega\zeta}{\zeta(\zeta+\pi p)},$$
and thus
\begin{equation}
\begin{array}{ccc}
{\displaystyle(\zeta I+\pi\Omega J)^{-1}(\zeta I-\pi\Omega J)=\left(\frac1\zeta I-\frac{\pi\Omega\zeta}{\zeta(\zeta+\pi p)}\right)(\zeta I-\pi\Omega J)}\\
{\displaystyle=I-\frac{\pi\Omega J}\zeta-\frac{\pi\Omega J}{\zeta+\pi p}+\frac{\pi^2p\Omega J}{\zeta(\zeta+\pi p)}=I-\frac{2\pi\Omega J}{\zeta+\pi p}.}
\end{array}\label{eq4.34}
\end{equation}
As a result, \eqref{eq4.34} implies the {\bf Riemann matrix boundary value problem} \cite{18,21}
\begin{equation}
\Phi_+(x)=B(x,\lambda)\Phi_-(x)+\frac{\Omega(x)J}{\zeta(x,\lambda)+\pi p(x)}\label{eq4.35}
\end{equation}
where
\begin{equation}
B(x,\lambda)=I-\frac{2\pi\Omega(x)J}{\zeta(x,\lambda)+\pi p(x)}.\label{eq4.36}
\end{equation}
By $X(z,\lambda)$, denote the {\bf canonical matrix function} which is the solution to the {\bf Riemann homogenous matrix problem}
\begin{equation}
X_+(x,\lambda)=B(x,\lambda)X_-(x,\lambda)\label{eq4.37}
\end{equation}
normalized by the condition $X(\infty,\lambda)=I$ where $X_\pm(x,\lambda)=X(x\pm i0,\lambda)$. Then \eqref{eq4.35} implies
\begin{equation}
X_+^{-1}(x,\lambda)\Phi_+(x)=X_-^{-1}(x,\lambda)\Phi_-(x)+\frac{X_+^{-1}(x,\lambda)\Omega(x)J}{\zeta(x,\lambda)+\pi p(x)}\label{eq4.38}
\end{equation}
Set the matrix function
\begin{equation}
\psi(z,\lambda)=\frac1{2\pi i}\int\limits_{\mathbb{R}}\frac{X_+^{-1}(x,\lambda)\Omega(x)J}{\zeta(x,\lambda)+\pi p(x)}\frac{dx}{x-z}\quad(z\in\mathbb{C}).\label{eq4.39}
\end{equation}
The functions $\omega_{\gamma,\gamma}$ \eqref{eq4.25} belong to $L^2(\mathbb{R})$ due to \eqref{eq4.30}, therefore Sokhotski formulas \eqref{eq4.31} are true for $\psi(z,\lambda)$ as a function of $z\in\mathbb{C}$ (with $\lambda$ fixed), therefore equation \eqref{eq4.38} becomes
$$X_+^{-1}(x,\lambda)\Phi_+(x)-\psi_+(x,\lambda)=X_-^{-1}(x,\lambda)\Phi_-(x)-\psi_-(x,\lambda)$$
where $\psi_\pm(x,\lambda)=\psi(x\pm x0,\lambda)$. So, the matrix function $X^{-1}(z,\lambda)\Phi(z)-\psi(z,\lambda)$ of $z$ can be holomorphically extended across $\mathbb{R}$ (for every fixed $\lambda$) and coincides with matrix-valued polynomial $P_\varkappa(z,\lambda)$ of $z$ \cite{18,21} of degree $\varkappa\geq0$ where $\varkappa$ is `index' \cite{18,20,21} of the Riemann problem \eqref{eq4.25}, therefore
\begin{equation}
\Phi(z)=\chi(z,\lambda)\{\psi(z,\lambda)+P_\varkappa(z,\lambda)\}.\label{eq4.40}
\end{equation}
Characteristic matrix function $S_\Delta(\lambda,B)$ is
$$S_\Delta(\lambda,B)=I-i[\langle\widetilde{\varphi}(\lambda I-\widetilde{X})f_\alpha,e_\beta\rangle]=I+i\left[\int\limits_\mathbb{R}\frac{f_\alpha(t)\overline{v_\beta}(t)}{t-\lambda}dt\right]$$
$$=I+2\pi\Phi(\lambda),$$
due to \eqref{eq4.24}, \eqref{eq4.29}, \eqref{eq4.27}. Taking into account \eqref{eq4.39}, \eqref{eq4.40}, one obtains
$$S_\Delta(\lambda,B)=I+i\int\limits_\mathbb{R}\chi(\lambda,\lambda)X_+^{-1}(x,\lambda)\frac{\Omega(x)J}{\zeta(x,\lambda)+\pi p(x)}\frac{dx}{x-\lambda}-2\pi\chi(\lambda,\lambda)P_\varkappa(\lambda,\lambda),$$
and since $S_\Delta(\infty,B)=I$ and $\chi(\infty,\lambda)=I$, then the last term in this sum vanishes.

\begin{theorem}\label{t4.2}
Let $X$, the spectral root of a quadratic pencil $L(\lambda)$ \eqref{eq3.2}, be a bounded self-adjoint operator with simple absolutely continuous spectrum ($\sigma(x)=[a,b]$, $-\infty<a<b<\infty$) and suppose that $\rank(\varphi^*\sigma\varphi)=r<\infty$ and $K(x,t)$ is the kernel in \eqref{eq3.15} for $\widetilde{\varphi^*\sigma\varphi}$ ($=\widetilde{K}$) corresponding to the spectral decomposition \eqref{eq3.21} of the operator $\widetilde{X}$ and $K(x,t)$ is given by \eqref{eq4.22} where $\{v_\alpha\}_1^r$ is a set of linearly independent functions from $L^2(a,b)$ for which \eqref{eq4.30} takes place and $j_{\alpha,\beta}=\langle Je_\beta,e_\alpha\rangle$ ($\sigma=J=J^*=J^{-1}$, $\{e_\alpha\}_1^r$ is the orthonormal basis in $E$).

Then the characteristic matrix function $S_\Delta(\lambda,B)=[\langle S_\Delta(\lambda,B)e_\alpha,e_\beta\rangle]$ corresponding to the model pair $\widetilde{X}$, $\widetilde{Y}$ \eqref{eq4.21} is
\begin{equation}
S_\Delta(\lambda,B)=I+i\int\limits_\mathbb{R}\chi(\lambda,\lambda)\chi_+^{-1}(x,\lambda)\frac{\Omega(x)J}{\zeta(x,\lambda)+\pi p(x)}\frac{dx}{x-\lambda},\label{eq4.41}
\end{equation}
here $\zeta(x,\lambda)=\lambda+x+b(x)$ ($b(x)$ is a real function and $\sup\limits_x|b(x)|<\infty$); $p(x)=\sum\limits_{\alpha,\beta=1}^rv_\alpha(x)j_{\alpha,\beta}\overline{v_\beta}(x)$; $\Omega(x)$ is given by \eqref{eq4.27}; $x(z,\lambda)$, for every fixed $\lambda$, is the canonical matrix function of problem \eqref{eq4.37} ($B(x,\lambda)$ is given by \eqref{eq4.36}) normalized by the condition $\chi(\infty,\lambda)=I$.
\end{theorem}
\vspace{5mm}

{\bf 4.3.} Give form of the characteristic function $S_\Delta(\lambda,B)$ obtained as a result of coupling of primaries.

Suppose that $\dim E=1$, $H=\mathbb{C}$, and
\begin{equation}
B=b_1\,(\in\mathbb{R}),\quad A=\lambda_1\,(\in\mathbb{C}_+),\quad\lambda_1-\overline{\lambda_1}=i\beta_1^2\,(\beta_1>0),\label{eq4.42}
\end{equation}
hence it follows that $\varphi=\beta_1$, $\sigma_1$, and characteristic function $S_{\Delta_1}(\lambda,B)$ is
\begin{equation}
S_{\Delta_1}(\lambda,B)=1-\frac{i\beta_1^2}{\lambda^2+\lambda b_1+\lambda_1}=\frac{\lambda^2+\lambda b_1+\overline{\lambda_1}}{\lambda^2+\lambda b_1+\lambda_1}.\label{eq4.43}
\end{equation}
Hence and from \eqref{eq1.27} it follows that
\begin{equation}
1-|S_{\Delta_1}(\lambda,B)|^2=\frac{\lambda-\overline{\lambda}}i(\lambda+\overline{\lambda}+b_1)\frac{\beta_1^2}{|\lambda^2+\lambda b_1+\lambda_1|^2},\label{eq4.44}
\end{equation}
therefore
\begin{equation}
1-|S_{\Delta_1}(\lambda,B)|^2=\left\{
\begin{array}{ccc}
>0,&{\displaystyle\Im\lambda\left(\Re\lambda+\frac{b_1}2\right)>0;}\\
=0,&{\displaystyle\Im\lambda=0\cup\Re\lambda=-\frac{b_1}2;}\\
<0,&{\displaystyle\Im\lambda\left(\Re\lambda+\frac{b_1}2\right)<0.}
\end{array}\right.\label{eq4.45}
\end{equation}
If
\begin{equation}
S_{\widetilde{\Delta}}(\lambda,\widetilde{B})=S_{\Delta_2}(\lambda,B_2)\cdot S(\lambda,B_1)\label{eq4.46}
\end{equation}
where
\begin{equation}
S_{\Delta_k}(\lambda,B_k)=\frac{\lambda^2+\lambda b_k+\overline{\lambda}_k}{\lambda^2+\lambda b_k+\lambda_k}\quad(b_k\in\mathbb{R},\lambda_k\in\mathbb{C}_+,k=1,2)\label{eq4.47}
\end{equation}
and, in accordance with \eqref{eq1.13}, \eqref{eq1.14},
\begin{equation}
\widetilde{B}=\left[
\begin{array}{ccc}
b_1&0\\
0&b_2
\end{array}\right],\quad\widetilde{A}=\left[
\begin{array}{ccc}
\lambda_1&0\\
i\beta_1\beta_2&\lambda_2
\end{array}\right].\label{eq4.48}
\end{equation}
Operator roots $\widetilde{X}$, $\widetilde{Y}$, due to \eqref{eq1.39}, of the pencil $\widetilde{L}(\lambda)=L_1(\lambda)L_2(\lambda)$ are expressed via roots $\{w_1^k,w_2^k\}$ of the equations $\lambda^2+\lambda b_k+\lambda_k=0$ ($k=1$, 2) by the formulas
\begin{equation}
\widetilde{X}=\left[
\begin{array}{ccc}
w_1^1&0\\
\gamma_{21}&w_1^2
\end{array}\right];\quad\widetilde{Y}=\left[
\begin{array}{ccc}
w_2^1&0\\
-\gamma_{21}&w_2^2
\end{array}\right],\label{eq4.49}
\end{equation}
besides, $w_1^k\in\mathbb{C}_+$, and $w_2^k\in\mathbb{C}_-$. Taking \eqref{eq1.40} for $\gamma_{21}$ into account, one obtains
\begin{equation}
w_2^2\gamma_{21}-\gamma_{21}w_1^1=i\beta_1\beta_2;\quad\gamma_{21}=i\frac{\beta_1\beta_2}{w_2^2-w_1^1},\label{eq4.50}
\end{equation}
and $w_2^2-w_1^1\not=0$ because $w_1^1\in\mathbb{C}_+$ and $w_2^2\in\mathbb{C}_-$. For $S_{\widetilde{\Delta}}(\lambda,\widetilde{B})$ \eqref{eq4.43}, for $b_1\geq b_2$, relations analogous to \eqref{eq4.45} are true,
\begin{equation}
1-|S_{\widetilde{\Delta}}(\lambda,\widetilde{B})|^2=\left\{
\begin{array}{lll}
>0;&{\displaystyle\left[\Im\lambda>0\cup\left(\Re\lambda+\frac{b_1}2\right)>0\right]\cup\left[\Im\lambda<0\cup\right.}\\
&{\displaystyle\left.\left(\Re\lambda+\frac{b_1}2\right)<0\right];}\\
=0;&{\displaystyle\lambda\in\left(-\infty,-\frac{b_1}2\right)\cup\left(-\frac{b_2}2,\infty\right);}\\
<0;&{\displaystyle\left[\Im\lambda>0\cup\left(\Re\lambda+\frac{b_1}2\right)<0\right]\cup\left[\Im\lambda<0\cup\right.}\\
&{\displaystyle\left.\left(\Re\lambda+\frac{b_1}2\right)>0\right]}
\end{array}\right.\label{eq4.51}
\end{equation}
Consider now the case when
\begin{equation}
S_{\widetilde{\Delta}}(\lambda,\widetilde{B})=\prod\limits_1^NS_{\Delta_k}(\lambda,B_k)\label{eq4.52}
\end{equation}
where $S_{\Delta_k}(\lambda,B_k)$ are given by \eqref{eq4.47} and $b_1\geq b_2\geq...\geq b_N$, then upon coupling one has
\begin{equation}
\widetilde{B}=\left[
\begin{array}{cccc}
b_1&0&...&0\\
0&b_2&...&0\\
...&...&...&...\\
0&0&...&b_N
\end{array}\right],\quad\widetilde{A}=\left[
\begin{array}{cccc}
\lambda_1&0&...&0\\
i\beta_1\beta_2&\lambda_2&...&0\\
...&...&...&...\\
i\beta_1\beta_N&i\beta_2\beta_N&...&\lambda_N
\end{array}\right].\label{eq4.53}
\end{equation}
Denote by $w_1^k$, $w_2^k$ the roots of equation $\lambda^2+\lambda b_k+\lambda_k=0$ ($1\leq k\leq N$) where $w_1^k\in\mathbb{C}_+$ and $w_2^k\in\mathbb{C}_-$. Equation \eqref{eq1.40} implies that the operators $\widetilde{X}$ and $\widetilde{Y}$ are
\begin{equation}
\widetilde{X}=\left[
\begin{array}{cccc}
w_1^1&0&...&0\\
\gamma_{21}&w_1^2&...&0\\
...&...&...&...\\
\gamma_{N,1}&\gamma_{N,2}&...&w_1^N
\end{array}\right];\quad\widetilde{Y}=\left[
\begin{array}{cccc}
w_2^1&0&...&0\\
-\gamma_{21}&w_2^2&...&0\\
...&...&...&...\\
-\gamma_{N,1}&-\gamma_{N,2}&...&w_2^N
\end{array}\right]\label{eq4.54}
\end{equation}
where the numbers $\gamma_{k,s}$ ($k>s$) are found (due to \eqref{eq1.40}) from the system of equations
\begin{equation}
(w_2^k-w_1^s)\gamma_{k,s}=i\beta_s\beta_k+\sum\limits_{s<l<k}\gamma_{k,l}\gamma_{l,s}.\label{eq4.55}
\end{equation}
The equation system \eqref{eq4.55} is solved step-by-step. First, the elements $\gamma_{k,k-1}$ situated below the main diagonal are found,
$$(w_2^k-w_1^{k-2})\gamma_{k,k-2}=i\beta_k\beta_{k-1}\quad(2\leq k\leq N),$$
knowing which one calculates the elements of the second diagonal below the main,
$$(w_2^k-w_1^{k-2})\gamma_{k,k-2}=i\beta_k\beta_{k-2}+\gamma_{k,k-1}\gamma_{k-1,k-2}\quad(3\leq k\leq N),$$
and finally
$$(w_2^N-w_1^1)\gamma_{N,1}=i\beta_N\beta_1+\sum\limits_{2<l<N}\gamma_{Nl}\cdot\gamma_{l,1}.$$
For $S_{\widetilde{\Delta}}(\lambda,\widetilde{B})$ \eqref{eq4.52}, analogously to \eqref{eq4.51},
\begin{equation}
1-|S_{\widetilde{\Delta}}(\lambda,\widetilde{B})|^2=\left\{
\begin{array}{lll}
>0;&{\displaystyle\left[\Im\lambda>0\cup\left(\Re\lambda+\frac{b_N}2\right)>0\right]\cup\left[\Im\lambda<0\cup\right.}\\
&{\displaystyle\left.\left(\Re\lambda+\frac{b_1}2\right)<0\right];}\\
=0;&{\displaystyle\lambda\in\left(-\infty,-\frac{b_1}2\right)\cup\left(-\frac{b_N}2,\infty\right);}\\
<0;&{\displaystyle\left[\Im\lambda>0\cup\left(\Re\lambda+\frac{b_1}2\right)<0\right]\cup\left[\Im\lambda<0\cup\right.}\\
&{\displaystyle\left.\left(\Re\lambda+\frac{b_N}2\right)>0\right].}
\end{array}
\right.\label{eq4.56}
\end{equation}
Equalities
\begin{equation}
S_{\Delta_k}(\lambda,B_k)=1-\frac{i\beta_k^2}{\lambda^2+\lambda b_k+\lambda_k}\quad(1\leq k\leq N)\label{eq4.57}
\end{equation}
imply that convergence of product \eqref{eq4.52} for $N\rightarrow\infty$ takes place if
\begin{equation}
\sum\limits_k^N\beta_k^2<\infty;\quad|b_k|<M<\infty\,(\forall k).\label{eq4.58}
\end{equation}
Model representations of $\widetilde{B}$, $\widetilde{A}$ and roots $\widetilde{X}$, $\widetilde{Y}$ become infinite matrices \eqref{eq4.53} and \eqref{eq4.54}, respectively, in the space $l^2(\mathbb{N})$. Boundedness of the operators $\widetilde{B}$, $\widetilde{A}$ and $\widetilde{X}$, $\widetilde{Y}$ follows from the conditions \eqref{eq4.58}.

Since $\lambda_k\rightarrow\mu_k\in\mathbb{R}$ ($k\rightarrow\infty$), due to \eqref{eq4.58}, then
$$1-\frac{i\beta_k^2}{\lambda^2+\lambda b_k+\mu_k}=\exp\left\{-\frac{i\beta_k^2}{\lambda^2+\lambda b_k+\mu_k}\right\}+O(\beta_k^4),$$
for limit points $\lambda_k\rightarrow\mu_k\in\mathbb{R}$, therefore the infinite product \eqref{eq4.52} corresponding to these $\mu_k$ is
$$\exp\left\{-i\int\limits_a^b\frac{d\sigma(t)}{\lambda^2+\lambda q(t)+\mu(t)}\right\}$$
and, after the substitution $\sigma(s)+c=t$, one obtains
\begin{equation}
S_\Delta(\lambda,B)=\exp\left\{-i\int\limits_0^l\frac{dt}{\lambda^2+\lambda b(t)+a(t)}\right\}\label{eq4.59}
\end{equation}
where $0\leq l<\infty$ and $b(t)$, $a(t)$ are real bounded functions.

Model representations of $\widetilde{A}$ and $\widetilde{B}$ corresponding to the characteristic function $S_\Delta(\lambda,B)$ in the space $L^2(0,l)$ are given by the formulas
\begin{equation}
(\widetilde{B}f)(x)=b(x)f(x);\quad(\widetilde{A}f)(x)=a(x)+i\int\limits_0^xf(t)dt.\label{eq4.60}
\end{equation}
The operator $\sigma=1$, and $\widetilde{\varphi}$ is
$$\widetilde{\varphi}:L^2(0,l)\rightarrow\mathbb{C};\quad\widetilde{\varphi}f=\int\limits_0^lf(x)dx\quad(f\in L^2(0,l),\,\varphi^*\xi=\xi).$$
Elementary calculations show that the characteristic function $S_{\widetilde{\Delta}}(\lambda,\widetilde{B})=I-i\widetilde{\varphi}\widetilde{L}(\lambda)\widetilde{\varphi}^*$ coincides with \eqref{eq4.59}.

Denote by $w_1(t)$, $w_2(t)$ roots of the equation
\begin{equation}
\lambda^2+\lambda b(t)+a(t)=0,\label{eq4.61}
\end{equation}
then
\begin{equation}
w_1(t)+w_2(t)=-b(t),\quad w_1(t)w_2(t)=a(t)\label{eq4.62}
\end{equation}
Roots $w_1(t)$, $w_2(t)$ may be either real or imaginary depending on the sign of determinant $D(t)=b^2(t)-4a(t)$. Define the operators
\begin{equation}
\begin{array}{ccc}
{\displaystyle(\widetilde{X}f)(x)=w_1(x)f(x)+i\int\limits_0^xK(x,t)f(t)dt;}\\
{\displaystyle(\widetilde{Y}f)(x)=w_2(x)f(x)-i\int\limits_0^xK(x,t)f(t)dt,}
\end{array}\label{eq4.63}
\end{equation}
then $\widetilde{X}+\widetilde{Y}=\widetilde{B}$ \eqref{eq4.60}, due to \eqref{eq4.63}. Condition $\widetilde{Y}\widetilde{X}=\widetilde{A}$ \eqref{eq4.60} implies
$$w_2(x)i\int\limits_0^xK(x,t)f(t)dt-i\int\limits_0^xK(x,t)w_1(t)f(t)dt$$
$$+\int\limits_0^xK(x,t)\int\limits_0^tK(t,s)f(s)dsdt=i\int\limits_0^xf(t)dt.$$
Since this equality should hold for all $f\in L^2(0,l)$, then one has equation for the kernel $K(x,t)$
\begin{equation}
[w_2(x)-w_1(t)]K(x,t)-i\int\limits_t^xK(x,\xi)K(\xi,t)d\xi=1.\label{eq4.64}
\end{equation}
Look for a solution to this equation in the form $K(x,t)=C(x)D(t)$, then
\begin{equation}
C(x)D(t)\left\{w_2(x)-w_1(t)-i\int\limits_t^xC(\xi)D(\xi)d\xi\right\}=1.\label{eq4.65}
\end{equation}
Hence, one finds that
\begin{equation}
C(x)D(x)=\frac1{w_2(x)-w_1(x)}.\label{eq4.66}
\end{equation}

\begin{remark}\label{r4.2}
If $w_2(x)=w_1(x)$, then \eqref{eq4.64} implies that
$$2w_1(x)=-b(x),\quad w_1(x)=a(x)^{1/2},$$
i. e., $a(x)^{1/2}=-b(x)/2$. Therefore, under condition that
\begin{equation}
2\sqrt{a(x)}\not=-b(x)\quad(\forall x),\label{eq4.67}
\end{equation}
expression \eqref{eq4.56} exists and is correct.
\end{remark}

Equation \eqref{eq4.61} implies
\begin{equation}
K(x,t)=\left\{w_2(x)-w_1(t)-i\int\limits_t^x\frac{d\xi}{w_2(\xi)-w_1(\zeta)}\right\}^{-1}.\label{eq4.68}
\end{equation}
So, if \eqref{eq4.67} takes place, then for $\widetilde{X}$, $\widetilde{Y}$ representations \eqref{eq4.63} are true, where $K(x,t)$ is given by \eqref{eq4.68}.

\begin{theorem}\label{t4.3}
Suppose that $\dim E=1$ and characteristic function $S_\Delta(\lambda,B)$ \eqref{eq1.26} is
\begin{equation}
S_\Delta(\lambda,B)=\prod\limits_1^N\frac{\lambda^2+\lambda b_k+\overline{\lambda}_k}{\lambda^2+\lambda b_k+\lambda_k}\cdot\exp\left\{-i\int\limits_0^l\frac{dt}{\lambda^2+\lambda b(t)+a(t)}\right\}\label{eq4.69}
\end{equation}
where $b_k\in\mathbb{R}$, $\lambda_k\in\mathbb{C}_+$ ($\lambda_k-\overline{\lambda}_k=i\beta_k^2$ and \eqref{eq4.53} takes place; $N\leq\infty$; $b(t)$, $a(t)$ are real bounded functions; $0\leq l<\infty$.

Then the operators $\widetilde{A}$, $\widetilde{B}$ \eqref{eq4.53} and $\widetilde{X}$, $\widetilde{Y}$ \eqref{eq4.54} are model representations of the infinite product in \eqref{eq4.52} and $\widetilde{A}$, $\widetilde{B}$ \eqref{eq4.60}, and $\widetilde{X}$, $\widetilde{Y}$ \eqref{eq4.63} are model representations of the second factor in \eqref{eq4.63}.
\end{theorem}

\section{Problems}

In the conclusion, we state problems which are left unsolved in this field.
\vspace{5mm}

{\bf 1.} If operator pairs $B$, $A$ in $H$ and $B'$, $A'$ in $H'$ are unitarily equivalent (i. e., there exists a unitary operator $U:$ $H\rightarrow H'$ such that $UA=A'U$ and $UB=B'U$), and thus the pencils $L(\lambda)$ and $L'(\lambda)$ are unitarily equivalent, then $S_\Delta(\lambda,B)=S_{\Delta'}(\lambda,B')$. Prove the inverse statement. If characteristic functions of {\bf simple} ?! pairs $B$, $A$ and $B'$, $A'$ coincide, $S_\Delta(\lambda,B)=S_{\Delta'}(\lambda,B')$, then the pairs $B$, $A$ and $B'$, $A'$ are unitarily equivalent. Notion of the simplicity of a pair $B$, $A$ requires explanation. When $B=0$ ($B'=0$), this statement is well-known \cite{23,24,35}.
\vspace{5mm}

{\bf 2.} Give description of the class of operator-valued functions $S(\lambda)$ which for some $B$ ($=B^*$) and $A$ ($A-A^*=i\varphi^*\sigma\varphi$) coincide with $S_\Delta(\lambda,B)$. This means characterization of the functions
$$S_\Delta(\lambda,B)=I-i\varphi L^{-1}(\lambda)\varphi^*\sigma.$$
For $B=0$, function $S_\Delta(\lambda,B)$ (after the substitution $\lambda^2=\lambda$) coincides with classical Livsic characteristic function, description of which is well-known (\cite{23,24,35}). This problem, after Caley transform, is equivalent to the description of a class of functions of the form $\varphi(\lambda I+\lambda B+A_R)\varphi^*$ where $B=B^*$, $A_R^*=A_R$ which generalizes class of Nevanlinna functions.

\begin{remark}\label{r5.1}
Problems 1, 2, probably, will require additional conditions on pairs $A$, $B$ (similar to commutation relations).
\end{remark}

{\bf 3.} Obtain, similar to Potapov $J$-theory, multiplicative decomposition for the functions $S_\Delta(\lambda,B)$.
\vspace{5mm}

{\bf 4.} Construct functional model of a pencil $L(\lambda)$ the way it was done \cite{24} for non-self-adjoint operators, i. e., for $B=0$. Here naturally arise problems of construction of dilation theory for pencils $L(\lambda)$.

\renewcommand{\refname}{ \rm 
\end{Large}
\end{document}